# Hilbert functions of graded algebras over Artinian rings[*]

Cristina Blancafort

Departament d'Àlgebra i Geometria

Universitat de Barcelona

Gran Via 585, 08007 Barcelona, Spain

e-mail: `blancafo@cerber.mat.ub.es`

March 1995

## Abstract

In this paper we give an effective characterization of Hilbert functions and polynomials of standard algebras over an Artinian equicharacteristic local ring; the cohomological properties of such algebras are also studied. We describe algorithms to check the admissibility of a given function or polynomial as a Hilbert function or polynomial, and to produce a standard algebra with a given Hilbert function.

## Introduction

Let $(R_0, \mathfrak{m}, k)$ be an Artinian local ring, $R = R_0[X_1, \ldots, X_b]$ and $I \subseteq R_+ = \oplus_{n \geq 1} R_n$ a homogeneous ideal. We will call standard $R_0$-algebra a graded algebra of the form $S = R/I$, and we will denote by $H_S(n) = \lambda_{R_0}(S_n)$ the Hilbert function of $S$. The study of Hilbert functions goes back a century in time; its origin is the celebrated result due to Hilbert, see [Hil90]:

[**Hilbert 1890**] If $R_0 = k$ is a field, then $H_S$ is asymptotically polynomical.

Later on Macaulay characterized Hilbert functions in the case $R_0 = k$ in [Mac27], [Mac16]:

[**Macaulay 1927**] Let $H : \mathbb{N} \to \mathbb{N}$ be a numerical function; then $H$ is the Hilbert function of a standard $k$-algebra if and only if $H(0) = 1$ and $H(n+1) \leq (H(n)_n)_+^+$ for all $n \geq 1$, see [BH93] for a proof. Afterwards, Samuel and Serre extended Hilbert's result to the Artinian case in [Sam51] and [Ser65]. In view of this situation, we found it natural to consider the following problems:

**(P1)** Extension of Macaulay's characterization to the Artinian case,

**(P2)** Characterization of the polynomials in $\mathbb{Q}[X]$ which are the Hilbert polynomial of a standard $R_0$-algebra.

---

Another interesting result in this line is Gotzmann's regularity theorem, see [Got78] and [Gre89]. This theorem, according to Green's presentation, provides us with an alternative expression of Hilbert polynomials more "combinatorical-like" than the usual one. It is deeply related to the study of the behaviour of Hilbert functions under hyperplane section, which is a standard method to perform inductive proofs in Commutative Algebra and Algebraic Geometry.

**[Gotzmann 1978, Green 1989]** If $R_0$ is a field, there exist uniquely determined integers $c_1 \geq c_2 \geq \cdots \geq c_s \geq 0$ such that the Hilbert polynomial of $S = k[X_1, \ldots, X_b]/I$ can be written

$$h_S(X) = \binom{X + c_1}{c_1} + \binom{X + c_2 - 1}{c_2} + \cdots + \binom{X + c_s - (s-1)}{c_s}.$$

Furthermore, the ideal sheaf $\mathcal{I}$ associated to $I$ is $s$-regular.

Hence an additional problem we have considered is:

**(P3)** Extension of Gotzmann's persistence theorem to the Artinian case and its relation with Castelnuovo-Mumford's regularity of the local cohomology $H^i_{S_+}(S)$.

The aim of this work is to study problems (P1), (P2) and (P3) in the case where $R_0$ is an Artinian $k$-algebra. Besides of being an interesting object of study in its own right, the theory of Hilbert functions of standard algebras over Artinian rings is the natural framework to study Hilbert functions of $\mathfrak{m}$-primary ideals in local rings, the Hilbert scheme and infinitesimal deformations. See for instance Remark 3.15 and Examples 5.6 and 5.7 for some results in this line. Further discussions on these fields will appear in a forthcoming paper.

In order to study the combinatorics of $R$ we introduce an ordered set of submodules of $R$ which considers both the combinatorics of the monomials and the structure of the base ring $R_0$. This set plays the role of a total ordering in the set of monomials of $k[X_1, \ldots, X_b]$ and extends the usual reverse lexicographical ordering. This will be especially neat when $R_0$ is a ring of deformations, i.e. $R_0 = k[\varepsilon]$.

The characterization theorem obtained in the study of problem (P1) takes into account the embedding dimension $b$ of the standard algebra $R_0$ such that $H = H_S$. Let us stress that this result is deeper than the mere generalization of Macaulay's theorem. The only information this straightforward generalization provides about the least possible value of $b$ is $b_{min} \leq H(1)$, which is enough when $R_0 = k$. In the general case we know only that $b_{min} \geq H(1)/\lambda_{R_0}(R_0)$, hence one needs to refine the conditions on $H$ to determine the minimal embedding dimension.

The extension of Gotzmann's result will provide us with bounds for the annihilation of the local cohomology $H^i_{S_+}(S)$. These bounds are computed in terms of the Hilbert-Samuel coefficients. For instance, we recover Hoa's result $a(S) \leq e(S) - \dim(S) - 1$, see Remark 3.13. We also compute the value $s$ appearing in Gotzmann's result and we show that it is a polynomial function of the Hilbert-Samuel coefficients, so we can make effective a result by Mumford about the regularity of ideal sheaves, see Remark 3.12.

We have been strongly concerned about the effectiveness of the results obtained. For instance: Macaulay's characterization as it is formulated is not an effective result, since there is no way to check the condition $H(n+1) \leq (H(n)_n)^+_+$ for all $n \in \mathbb{N}$. We will describe an algorithm to check these conditions in a finite number of steps, for any asymptotically polynomial function $H$. Specifically, we give algorithms to determine:

(i) whether a polynomial $P \in \mathbb{Q}[X]$ is a Hilbert polynomial,

(ii) whether an asymptotically polynomical function is a Hilbert function,

(iii) the minimal embedding dimension for a realizing algebra of a Hilbert function,

and we also compute a generating system, that will be minimal in the case $R_0 = k$, for the ideal $I$ such that a realizing algebra is $S = R/I$. In the general case the generating system will depend on a compositon series on $R_0$. Nevertheless, in the case where $R_0$ is finitely generated as a $k$-algebra, we can always obtain a composition series via Gröbner basis, see [CLO92], Proposition 1(ii) in Chapter 5, §3.

Let us describe the organization of the paper: In Section 1 we fix some notations we will use throughout the subsequent sections. In Section 2 we give a characterization of Hilbert functions of standard $R_0$-algebras for an equicharacteristic ring $R_0$, fixing the embedding dimension of the standard algebra $S$ such that $H = H_S$. The proof of the statements in this section and in the following one is inspired in the ideas exposed by M. Green in [Gre89]. In section 3 we prove an improved extension of Gotzmann's regularity theorem to the Artinian equicharacteristic case. This theorem provides us with an expression of Hilbert polynomials which will be useful to characterize them and to make our results effective, see Sections 4 and 5. We also obtain bounds for the annihilation of the local cohomology. In Section 4 we give an effective characterization of Hilbert polynomials of standard $R_0$-algebras with $R_0$ equicharacteristic. Again we formulate a version with a given embedding dimension which will be the key tool for the last section. Finally, in Section 6 we collect all the preceding statements in order to describe the algorithms which make them effective.

**Acknowledgements.** The author wishes to thank J. Elias for introducing her to the theory of Hilbert functions and for many useful discussions about this material. The members of the Seminari d'Àlgebra Commutativa de Barcelona are also to be thanked.

## 1 Notations

Let $R = \oplus_{n \geq 0} R_n$ be a $d$-dimensional graded ring such that $R_0$ is an Artinian local ring and $R$ is a $R_0$-module finitely generated by $R_1$; from now on we will call such a ring a standard $R_0$-algebra. As usual, we will denote by $R_+ = \oplus_{n \geq 1} R_n$ the irrelevant ideal of $R$. Let $H_R(n) := \lambda_{R_0}(R_n)$ be the Hilbert function of $R$. It is well known that, for $n \gg 0$, $H_R$ coincides with a polynomial $h_R$ of degree $d-1$ which is called the Hilbert polynomial of $R$. This and the fact that $H_R(n) \in \mathbb{N}$ for all $n$ suggest the definitions

$$\mathbb{Q}[X; \mathbb{Z}] = \{P \in \mathbb{Q}[X] \mid P(n) \in \mathbb{Z} \text{ for } n \gg 0\}$$

$$\mathbb{Q}[X; \mathbb{N}] = \{P \in \mathbb{Q}[X] \mid P(n) \in \mathbb{N} \text{ for } n \gg 0\}.$$

Clearly $\mathbb{Q}[X; \mathbb{N}]$ is the set of polynomials in $\mathbb{Q}[X; \mathbb{Z}]$ which have a positive leading coefficient. It is also known that $\mathbb{Q}[X; \mathbb{Z}]$ is a free $\mathbb{Z}$-module with basis $\{\binom{X+i}{i} \mid i \geq 0\}$, where $\binom{X+i}{i} = (X+i)(X+i-1)\ldots(X+1)/i!$ and $\binom{X}{0} = 1$. As a consequence, we can express

$$h_R(X) = e_0 \binom{X+d-1}{d-1} - e_1 \binom{X+d-2}{d-2} + \cdots + (-1)^{d-2} e_{d-2} \binom{X+1}{1} + (-1)^{d-1} e_{d-1}$$

with $e_i = e_i(R) \in \mathbb{Z}$ and $e_0 > 0$. The $e_i$ are called the normalized Hilbert-Samuel coefficients of $R$. We define also the regularity index of $R$ as

$$i(R) = \min\{k \in \mathbb{N} \mid H_R(n) = h_R(n) \ \forall \ n \geq k\}.$$

Similarly, if $(A, \mathfrak{m})$ is a $d$-dimensional local ring and $\mathfrak{a}$ is a $\mathfrak{m}$-primary ideal, we define $H^0_\mathfrak{a}(n) := \lambda_A(\mathfrak{a}^n/\mathfrak{a}^{n+1})$ and call it the Hilbert function of $\mathfrak{a}$; it coincides, for $n \gg 0$, with a polynomial $h^0_\mathfrak{a}$ of degree $d-1$ which we will call the Hilbert polynomial of $\mathfrak{a}$. In the same way as before we have that $h^0_\mathfrak{a} \in \mathbb{Q}[X;\mathbb{N}]$, and so we can write

$$h^0_\mathfrak{a}(X) = e_0 \binom{X+d-1}{d-1} - e_1 \binom{X+d-2}{d-2} + \cdots + (-1)^{d-2} e_{d-2} \binom{X+1}{1} + (-1)^{d-1} e_{d-1}$$

with $e_i = e_i(\mathfrak{a}) \in \mathbb{Z}$ and $e_0 > 0$. The $e_i$ are called the normalized Hilbert-Samuel coefficients of $\mathfrak{a}$. We define also the regularity index of $H^0_\mathfrak{a}$ as

$$i(H^0_\mathfrak{a}) = \min\{k \in \mathbb{N} \mid H^0_\mathfrak{a}(n) = h^0_\mathfrak{a}(n) \ \forall \ n \geq k\}.$$

In the case $\mathfrak{a} = \mathfrak{m}$ we will write $H^0_A$, $h^0_A$ and $e_i(A)$ instead of $H^0_\mathfrak{m}$, $h^0_\mathfrak{m}$ and $e_i(\mathfrak{m})$.

The relationship between the two cases is the following: let $gr_\mathfrak{a}(A) = \oplus_{n \geq 0} \mathfrak{a}^n/\mathfrak{a}^{n+1}$ denote the associated graded ring of $\mathfrak{a}$. Then $gr_\mathfrak{a}(A)_0 = A/\mathfrak{a}$ is an Artinian local ring and $gr_\mathfrak{a}(A)$ is a standard $(A/\mathfrak{a})$-algebra. In this way we have $H^0_\mathfrak{a} = H_{gr_\mathfrak{a}(A)}$. Reciprocally, if we have a standard $R_0$-algebra $R$ and $\mathfrak{M}$ is the maximal homogeneous ideal of $R$, the local ring $A = R_\mathfrak{M}$ and the $\mathfrak{M}A$-primary ideal $\mathfrak{a} = R_+ A$ verify $gr_\mathfrak{a}(A) \cong R$ as graded rings. Hence by the previous remark $H_R = H^0_\mathfrak{a}$.

## 2 Characterization of Hilbert functions

The main goal of this section is to obtain a characterization theorem for Hilbert functions of graded algebras over Artinian equicharacteristic rings, see Theorem 2.8. We remark that our result is stronger than the straightforward generalization of the classical Macaulay's result. This straightforward generalization is obtained as Corollary 2.10 and does not determine the embedding dimension of the realizing algebra when the base ring is not a field. This is why the conditions which characterize Hilbert functions of algebras *with a given embedding dimension* need to be refined.

We need some combinatorics to proceed: Given integers $n, d \geq 1$, it is known that there exist uniquely determined integers $k_d > k_{d-1} > \ldots > k_\delta \geq \delta \geq 1$ such that

$$n = \binom{k_d}{d} + \binom{k_{d-1}}{d-1} + \cdots + \binom{k_\delta}{\delta}.$$

This is called the $d$-binomial expansion of $n$. We define then

$$(n_d)^+ = \binom{k_d+1}{d} + \binom{k_{d-1}+1}{d-1} + \cdots + \binom{k_\delta+1}{\delta}, \qquad (n_d)^- = \binom{k_d-1}{d} + \binom{k_{d-1}-1}{d-1} + \cdots + \binom{k_\delta-1}{\delta},$$

$$(n_d)^+_+ = \binom{k_d+1}{d+1} + \binom{k_{d-1}+1}{d} + \cdots + \binom{k_\delta+1}{\delta+1}, \qquad (n_d)^-_- = \binom{k_d-1}{d-1} + \binom{k_{d-1}-1}{d-2} + \cdots + \binom{k_\delta-1}{\delta-1},$$

with the convention that $\binom{i}{j} = 0$ if $i < j$ and $\binom{i}{0} = 1$ for all $i \geq 0$. We also define $(0_d)^+ = (0_d)^- = (0_d)^+_+ = (0_d)^-_- = 0$ for all $d \geq 1$. Notice that we immediately obtain the $d$ and $(d+1)$-binomial expansions of $(n_d)^+$ and $(n_d)^+_+$ respectively. We will use as simplified notation $((n_d)^-)^+_+$ instead of $(((n_d)^-)_d)^+_+$.

We refer the reader to [Rob90] for some properties of these functions that will be used in the sequel. However, for the reader's convenience we list here the most used ones:

**Lemma 2.1** Let $n = \binom{k_d}{d} + \binom{k_{d-1}}{d-1} + \cdots + \binom{k_\delta}{\delta}$ and $m = \binom{l_d}{d} + \binom{l_{d-1}}{d-1} + \cdots + \binom{l_\mu}{\mu}$ be the d-binomial expansions of $n, m \geq 1$. Then we have:

(i) Define $k_{\delta-1} = \cdots = k_1 = 0$, $l_{\mu-1} = \cdots = l_1 = 0$. Then $n \leq m$ if and only if $(k_d, \ldots, k_1) \leq (l_d, \ldots, l_1)$ in the lexicographical order.

(ii) If $k_\delta > \delta$, then it holds that $((n-1)_d)^- < (n_d)^-$.

(iii) If $n < m$ then $(n_d)^+ < (m_d)^+$, $(n_d)^- \leq (m_d)^-$, $(n_d)^+_+ < (m_d)^+_+$ and $(n_d)^-_- \leq (m_d)^-_-$

(i) and (iii) can be found in [Rob90], §4, while (ii) is Lemma 4.2.11 (b) in [BH93].

In order to prove the characterization theorem, the following result will assure us of the existence of "good" linear forms that will allow us to perform the inductive step. The linear forms described in (i) are the best possible in the following sense: if $V = I_d$, $I$ a homogeneous ideal, they are the elements in $R_1 \cap k^b$ which are closest to being non-zero divisors in $R/I$. This will be made clear in the proof of Lemma 3.3.

**Theorem 2.2** Let $(R_0, \mathfrak{m})$ be an Artinian local equicharacteristic ring with infinite residue field $k$ and $R = R_0[X_1, \ldots, X_b]$. Consider an integer $d \geq 1$ and a $R_0$-submodule $V \subseteq R_d$. Then

(i) $U_R(d, V) = \{h \in R_1 \cap k^b \mid \lambda_{R_0}(V + hR_{d-1}) \text{ is maximal}\}$ is a Zariski-open in $R_1 \cap k^b$.

(ii) Let $\lambda_{R_0}(R_d/V) = \binom{d+b-1}{b-1}q + r$ be the Euclidean division of $\lambda_{R_0}(R_d/V)$ by $\binom{d+b-1}{b-1}$. For $h \in U_R(d, V)$, we put $\overline{R} = R/(h)$, $\overline{V} = (V + hR_{d-1})/hR_{d-1}$. Then we have

$$\lambda_{R_0}(\overline{R}_d/\overline{V}) \leq \binom{d+b-2}{b-2}q + (r_d)^-.$$

**Proof.** We first note that, being $R_0$ a complete local equicharacteristic ring, it contains a coefficient field which we will denote also by $k$, see for example [Mat86], Theorem 28.3.

Denote by $R_1 \cap k^b$ the set of linear forms in $R$ having all their coefficients in $k$. It can be naturally identified with $k^b$, and we will consider it as a topological space endowed with the Zariski topology.

Notice that for any finitely generated $R_0$-module $M$ the inclusion $k \subseteq R_0$ induces a $k$-vector space structure in $M$ and $\dim_k(M) = \lambda_{R_0}(M)$; let then $m$ be the maximal value of $\dim_k((V + hR_{d-1})/V)$ when $h \in R_1 \cap k^b$. Then $U_R(d, V)$ is the set of linear forms $h \in R_1 \cap k^b$ such that $\dim_k((V + hR_{d-1})/V) = m$.

Now consider for $h \in R_1 \cap k^b$ the $k$-linear map multiplication by $h = a_1X_1 + \ldots + a_bX_b$

$$\cdot h : R_{d-1} \to R_d/V.$$

Given $k$-bases of $R_{d-1}$ and $R_d/V$, we can describe this map by a matrix $M$ which entries are polynomical functions on $a_1, \ldots, a_b$. Since the image of $\cdot h$ is $(V + hR_{d-1})/V$, the complement of $U_R(d, V)$ in $k^b$ is the variety given by the ideal of $m \times m$ minors of $M$, $I_m(M)$.

Let $s = \lambda_{R_0}(R_0)$. To prove (ii) we will proceed by induction on $(b,d)$ in the lexicographical order.

Assume $b = 1$, then $R = R_0[X_1]$ and for all $h \in U_R(d,V)$ we have $h = a_1 X_1$ with $a_1 \in k$. Hence $\overline{R} = R/(h) = R_0$ and $\lambda_{R_0}(\overline{R}_d/\overline{V}) = 0$ since $d \geq 1$.

Assume now $d = 1$, in this case $V \subseteq R_1 = R_0 \langle X_1, \ldots, X_b \rangle$ and $\lambda_{R_0}(R_1/V) = bq + r$ with $0 \leq r < b$; notice that $q \leq s$. Since $h \notin \mathfrak{m}[X_1, \ldots, X_b]$, the multiplication by $h$ induces an isomorphism $R_0 \cong R_0 h$; from this and the isomorphism $\frac{V+R_0h}{V} \cong \frac{R_0h}{V \cap R_0h}$ we deduce that $U_R(1,V)$ coincides with the set of elements $h \in R_1 \cap k^b$ such that $\lambda_{R_0}(V \cap R_0 h)$ is minimal.

We will distinguish two cases:
(1) There exists $h_0 \in R_1 \cap k^b$ such that $\lambda_{R_0}(V \cap R_0 h_0) < s - q$. Then for all $h \in U_R(1,V)$ we must have $\lambda_{R_0}(V \cap R_0 h) < s - q$, and therefore

$$\begin{aligned}
\lambda_{R_0}(\overline{V}) &= \lambda_{R_0}((V + R_0 h)/R_0 h) \\
&= \lambda_{R_0}(V/(V \cap R_0 h)) = \lambda_{R_0}(V) - \lambda_{R_0}(V \cap R_0 h) \\
&> \lambda_{R_0}(V) - (s - q) \\
&= \lambda_{R_0}(R_1) - \lambda_{R_0}(R_1/V) - s + q \\
&= bs - (bq + r) - s + q = (b-1)s - (b-1)q - r.
\end{aligned}$$

Then $\lambda_{R_0}(\overline{V}) \geq (b-1)s - (b-1)q - r + 1$, hence
$\lambda_{R_0}(\overline{R}_1/\overline{V}) = \lambda_{R_0}(\overline{R}_1) - \lambda_{R_0}(\overline{V}) \leq (b-1)s - (b-1)s + (b-1)q + r - 1 = (b-1)q + (r_1)^-$.

(2) For all $h \in R_1 \cap k^b$ we have $\lambda_{R_0}(V \cap R_0 h) \geq s - q$. In particular $\lambda_{R_0}(V \cap R_0 X_i) \geq s - q$ for all $1 \leq i \leq b$, and since $\oplus_{i=1}^{b}(V \cap R_0 X_i) \subseteq V$ we deduce

$$\begin{aligned}
bq + r &= \lambda_{R_0}(R_1/V) \\
&\leq \lambda_{R_0}(R_1/(\oplus_{i=1}^{b}(V \cap R_0 X_i))) = \sum_{i=1}^{b} \lambda_{R_0}(R_0 X_i/(V \cap R_0 X_i)) \\
&= \sum_{i=1}^{b}(\lambda_{R_0}(R_0 X_i) - \lambda_{R_0}(V \cap R_0 X_i)) \\
&\leq bs - b(s-q) = bq
\end{aligned}$$

So, we get $r = 0$ and all the inequalities must be equalities. In particular, $\lambda_{R_0}(V \cap R_0 X_i) = s - q$ for all $1 \leq i \leq b$, and since it is by hypothesis the least possible value, $U_R(1,V)$ must be the set of linear forms $h \in R_1 \cap k^b$ such that $\lambda_{R_0}(V \cap R_0 h) = s - q$.
Let $h$ be a linear form belonging to $U_R(1,V)$, then it holds

$$\begin{aligned}
\lambda_{R_0}(\overline{R}_1/\overline{V}) &= \lambda_{R_0}(\overline{R}_1) - \lambda_{R_0}(\overline{V}) = (b-1)s - \lambda_{R_0}((V+R_0h)/R_0h) \\
&= (b-1)s - \lambda_{R_0}(V/(V \cap R_0 h)) = (b-1)s - \lambda_{R_0}(V) + \lambda_{R_0}(V \cap R_0 h) \\
&= (b-1)s - \lambda_{R_0}(R_1) + \lambda_{R_0}(R_1/V) + s - q \\
&= (b-1)s - bs + bq + s - q = (b-1)q = (b-1)q + (0_1)^-.
\end{aligned}$$

Finally, assume $b, d \geq 2$. Let $h \in U_R(d, V)$: we will denote by an overline the equivalence modulo $h$ and by $\pi : R \to \overline{R}$ the projection. Notice that, after a change of variables, we can consider $\overline{R}$ as a polynomial ring in $b-1$ variables. Then we have $\pi(R_1 \cap k^b) = \overline{R}_1 \cap k^{b-1}$.

Let us define $(V : h) = \{f \in R_{d-1} \mid hf \in V\}$ and consider the Zariski-open subset of $R_1 \cap k^b$

$$B = U_R(d-1, (V : h)) \cap \pi^{-1}(U_{\overline{R}}(d, \overline{V})).$$

Notice that these are Zariski-open with $k$ an infinite field, so $B$ must be nonempty. Pick $l \in B$, and denote by a hat accent the equivalence classes modulo $l$. If we define $((V : h) : l) = \{f \in R_{d-2} \mid lf \in (V : h)\}$ we have $((V : h) : l) = ((V : l) : h)$.

Let us list the following exact sequences that we will consider later on:

$$0 \to R_{d-1}/(V : h) \xrightarrow{\cdot h} R_d/V \to \overline{R}_d/\overline{V} \to 0 \tag{1}$$

$$0 \to R_{d-1}/(V : l) \xrightarrow{\cdot l} R_d/V \to \widehat{R}_d/\widehat{V} \to 0 \tag{2}$$

$$0 \to R_{d-2}/((V : l) : h) \xrightarrow{\cdot h} R_{d-1}/(V : l) \to \overline{R}_{d-1}/\overline{(V : l)} \to 0 \tag{3}$$

$$0 \to R_{d-2}/((V : h) : l) \xrightarrow{\cdot l} R_{d-1}/(V : h) \to \widehat{R}_{d-1}/\widehat{(V : h)} \to 0 \tag{4}$$

$$\overline{R}_{d-1}/\overline{(V : l)} \xrightarrow{\cdot \overline{l}} \overline{R}_d/\overline{V} \to \widehat{\overline{R}}_d/\widehat{\overline{V}} \to 0 \tag{5}$$

(5) being obtained from (2) modulo $h$.

Since $h \in U_R(d, V)$ and $l \in R_1 \cap k^b$, we have $\lambda_{R_0}(V + lR_{d-1}) \leq \lambda_{R_0}(V + hR_{d-1})$, and then from (1) and (2) we obtain $\lambda_{R_0}(R_{d-1}/(V : l)) \leq \lambda_{R_0}(R_{d-1}/(V : h))$. Applying to (3) and (4) together with the fact that $((V : l) : h) = ((V : h) : l)$ we get

$$\lambda_{R_0}(\overline{R}_{d-1}/\overline{(V : l)}) \leq \lambda_{R_0}(\widehat{R}_{d-1}/\widehat{(V : h)}).$$

Let us consider the following Euclidean divisions:

$$\lambda_{R_0}(R_d/V) = \binom{d+b-1}{b-1}q + r \quad \text{with } q \geq 0 \text{ and } 0 \leq r < \binom{d+b-1}{b-1},$$

$$\lambda_{R_0}(\overline{R}_d/\overline{V}) = \binom{d+b-2}{b-2}\overline{q} + \overline{r} \quad \text{with } \overline{q} \geq 0 \text{ and } 0 \leq \overline{r} < \binom{d+b-2}{b-2},$$

$$\lambda_{R_0}(R_{d-1}/(V : h)) = \binom{d-1+b-1}{b-1}\tilde{q} + \tilde{r} \quad \text{with } \tilde{q} \geq 0 \text{ and } 0 \leq \tilde{r} < \binom{d-1+b-1}{b-1}.$$

Recall that $b \geq 2$, so from Lemma 2.1 (ii) we have

$$0 \leq (r_d)^- \leq \left(\left(\binom{d+b-1}{d} - 1\right)_d\right)^-$$

$$< \left(\binom{d+b-1}{d}_d\right)^-$$

$$= \binom{d+b-2}{d} = \binom{d+b-2}{b-2}.$$

So the expression $\binom{d+b-2}{b-2}q + (r_d)^-$ is an Euclidean division, and hence the inequality

$$\binom{d+b-2}{b-2}\overline{q} + \overline{r} \leq \binom{d+b-2}{b-2}q + (r_d)^-,$$

which is what we want to prove, is equivalent to $(\overline{q}, \overline{r}) \leq (q, (r_d)^-)$ in the lexicographical order, since both sides are Euclidean divisions.

By (5) we have
$$\lambda_{R_0}(\overline{R}_d/\overline{V}) \le \lambda_{R_0}(\overline{R}_{d-1}/\overline{(V:l)}) + \lambda_{R_0}(\widehat{\overline{R}}_d/\widehat{\overline{V}})$$
$$\le \lambda_{R_0}(\widehat{R}_{d-1}/\widehat{(V:h)}) + \lambda_{R_0}(\widehat{\overline{R}}_d/\widehat{\overline{V}}).$$

Since $l \in U_R(d-1,(V:h))$ we can apply the induction hypothesis (on $d$) to the first term. Since $\bar{l} \in U_{\overline{R}}(d,\overline{V})$ we can apply the induction hypothesis (on $b$) to the second one. Therefore
$$\lambda_{R_0}(\overline{R}_d/\overline{V}) \le \binom{d-1+b-2}{b-2}\tilde{q} + (\tilde{r}_{d-1})^- + \binom{d+b-3}{b-3}\overline{q} + (\overline{r}_d)^-,$$

that is
$$\binom{d+b-2}{b-2}\overline{q} + \overline{r} \le \binom{d-1+b-2}{b-2}\tilde{q} + (\tilde{r}_{d-1})^- + \binom{d+b-3}{b-3}\overline{q} + (\overline{r}_d)^-.$$

From [Rob90], Corollary 4.6(a), we have $\overline{r} - (\overline{r}_d)^- = (\overline{r}_d)_-^-$, so the last inequality is equivalent to
$$\binom{d+b-3}{b-2}\overline{q} + (\overline{r}_d)_-^- \le \binom{d-1+b-2}{b-2}\tilde{q} + (\tilde{r}_{d-1})^- \quad (I)$$

As before $0 \le (\tilde{r}_{d-1})^- < \binom{d-1+b-2}{b-2}$ and the right-hand side in (I) is an Euclidean division. We need to distinguish two cases, which will correspond respectively to the case when the left-hand side is also an Euclidean division and the opposite one:

**Case (1):** $\overline{r} < \binom{d+b-2}{b-2} - b + 2$. Then $\overline{r} \le \binom{d+b-2}{d} - b + 1 = \binom{d+b-3}{d} + \binom{d+b-4}{d-1} + \cdots + \binom{b-1}{2}$ and this is a $d$-binomial expansion, so
$$(\overline{r}_d)_-^- \le \binom{d+b-4}{d-1} + \binom{d+b-5}{d-2} + \cdots + \binom{b-2}{1} = \binom{d+b-3}{d-1} - 1 < \binom{d+b-3}{b-2}$$

and then both sides of the inequality (I) are Euclidean divisions; therefore $(\overline{q},(\overline{r}_d)_-^-) \le (\tilde{q},(\tilde{r}_{d-1})^-)$ in the lexicographical order.

From the exact sequence (1) we obtain
$$\binom{d+b-1}{b-1}q + r = \binom{d+b-2}{b-2}\overline{q} + \overline{r} + \binom{d-1+b-1}{b-1}\tilde{q} + \tilde{r}$$
$$\ge \binom{d+b-2}{b-2}\overline{q} + \binom{d-1+b-1}{b-1}\overline{q} + \tilde{r} + \overline{r},$$

this is,
$$\binom{d+b-1}{b-1}q + r \ge \binom{d+b-1}{b-1}\overline{q} + \tilde{r} + \overline{r}.$$

Since $0 \le \tilde{r} + \overline{r} < \binom{d+b-2}{b-2} + \binom{d-1+b-1}{b-1} = \binom{d+b-1}{b-1}$, both sides of the inequality are Euclidean divisions, so $(\overline{q}, \tilde{r} + \overline{r}) \le (q, r)$ in the lexicographical order.

If $\overline{q} < q$ then case (1) is complete; assume then $\overline{q} = q$ (and therefore $\tilde{r} + \overline{r} \le r$). We must prove in this case $\overline{r} \le (r_d)^-$. Recall that $(\overline{q},(\overline{r}_d)_-^-) \le (\tilde{q},(\tilde{r}_{d-1})^-)$ in the lexicographical order and let us distinguish two subcases:

*Subcase (a):* $\overline{q} = \tilde{q}$. In this case we have $(\overline{r}_d)_-^- \le (\tilde{r}_{d-1})^- \le ((r-\overline{r})_{d-1})^-$, and by [Rob90], Corollary 4.11, we get $(\overline{r}_d)_- \le r - \overline{r}$. Therefore $(\overline{r}_d)^+ = \overline{r} + (\overline{r}_d)_- \le r$, and from [Rob90], Corollary 4.12, this inequality is equivalent to the desired one.

*Subcase (b)*: $\bar{q} < \tilde{q}$. Assume $r \geq \binom{d+b-1}{b-1} - b = \binom{d+b-2}{d} + \binom{d+b-3}{d-1} + \cdots + \binom{b}{2}$, then

$$(r_d)^- \geq \binom{d+b-3}{d} + \binom{d+b-4}{d-1} + \cdots + \binom{b-1}{2} = \binom{d+b-2}{b-2} - b + 1 \geq \bar{r}$$

and we immediately get $\bar{r} \leq (r_d)^-$. Therefore, we may assume $r = \binom{d+b-1}{b-1} - i$ with $b+1 \leq i \leq \binom{d+b-1}{b-1}$. From the exact sequence (1) we obtain

$$\binom{d+b-1}{b-1} q + \binom{d+b-1}{b-1} - i = \binom{d+b-2}{b-2} \bar{q} + \bar{r} + \binom{d-1+b-1}{b-1} \tilde{q} + \tilde{r}$$

$$\geq \binom{d+b-2}{b-2} \bar{q} + \bar{r} + \binom{d-1+b-1}{b-1}(\bar{q}+1) + \tilde{r}$$

$$= \binom{d+b-1}{b-1} \bar{q} + \binom{d-1+b-1}{b-1} + \bar{r} + \tilde{r}.$$

Since $q = \bar{q}$ we have

$$\binom{d+b-1}{b-1} - i \geq \binom{d-1+b-1}{b-1} + \bar{r} + \tilde{r},$$

hence

$$\bar{r} \leq \binom{d+b-1}{b-1} - \binom{d-1+b-1}{b-1} - i - \tilde{r} \leq \binom{d-1+b-1}{b-2} - i.$$

Since $r < \binom{d+b-1}{d}$ we have $(r_d)^-_- \leq \binom{d+b-2}{d-1}$, that is, $r - (r_d)^- \leq \binom{d+b-2}{b-1}$. Hence

$$\binom{d+b-1}{b-1} - i - (r_d)^- = r - (r_d)^- \leq \binom{d+b-2}{b-1}$$

i.e.

$$\bar{r} \leq \binom{d+b-2}{b-2} - i \leq (r_d)^-$$

and case (1) is complete.

**Case (2):** $\bar{r} \geq \binom{d+b-2}{b-2} - b + 2$. Then

$$\bar{r} \geq \binom{d+b-2}{d} - b + 2 = \binom{d+b-3}{d} + \binom{d+b-4}{d-1} + \cdots + \binom{b-1}{2} + \binom{1}{1}$$

and since this is a $d$-binomial expansion we have

$$(\bar{r}_d)^-_- \geq \binom{d+b-4}{d-1} + \binom{d+b-5}{d-2} + \cdots + \binom{b-2}{1} + \binom{0}{0} = \binom{d+b-3}{d-1}.$$

Therefore from (I) we obtain

$$\binom{d+b-3}{b-2}(\bar{q}+1) \leq \binom{d-1+b-2}{b-2} \tilde{q} + (\tilde{r}_{d-1})^-.$$

Notice that both sides are Euclidean divisions, so $\bar{q} + 1 \leq \tilde{q}$.

Let $\overline{r} = \binom{d+b-2}{b-2} - i$ with $1 \leq i \leq b-2$. From the exact sequence (1) we obtain

$$\binom{d+b-1}{b-1}q + r = \binom{d+b-2}{b-2}\overline{q} + \overline{r} + \binom{d-1+b-1}{b-1}\tilde{q} + \tilde{r}$$

$$= \binom{d+b-2}{b-2}(\overline{q}+1) - i + \binom{d-1+b-1}{b-1}\tilde{q} + \tilde{r}$$

$$\geq \binom{d+b-2}{b-2}(\overline{q}+1) - i + \binom{d-1+b-1}{b-1}(\overline{q}+1) + \tilde{r}$$

$$= \binom{d+b-1}{b-1}(\overline{q}+1) + \tilde{r} - i,$$

this is

$$\binom{d+b-1}{b-1}q + r + i \geq \binom{d+b-1}{b-1}(\overline{q}+1) + \tilde{r},$$

and since $\binom{d-1+b-1}{b-1} < \binom{d+b-1}{b-1}$ the right-hand side of this inequality is an Euclidean division. If $r < \binom{d+b-1}{b-1} - i$ then the left-hand side is also an Euclidean division and therefore $\overline{q} < q$ as we wanted to prove. Assume $r \geq \binom{d+b-1}{b-1} - i$, then since $r < \binom{d+b-1}{b-1}$ and $1 \leq i \leq b - 2 < \binom{d+b-1}{b-1}$ we get

$$\binom{d+b-1}{b-1}(\overline{q}+1) + \tilde{r} \leq \binom{d+b-1}{b-1}q + r + i < \binom{d+b-1}{b-1}(q+1) + i$$

and now both sides are Euclidean divisions, so $\overline{q} \leq q$. Then, it only remains to show that if $\overline{q} = q$ then $\overline{r} \leq (r_d)^-$.

Since $1 \leq i \leq b-2$, in the inequality $r \geq \binom{d+b-1}{d} - i = \binom{d+b-2}{d} + \binom{d+b-3}{d-1} + \cdots + \binom{b}{2} + \binom{b-i}{1}$, the last expression is a $d$-binomial expansion. Hence

$$(r_d)^- \geq \binom{d+b-3}{d} + \binom{d+b-4}{d-1} + \cdots + \binom{b-1}{2} + \binom{b-i-1}{1} = \binom{d+b-2}{b-2} - i = \overline{r}$$

and this completes the proof. $\square$

As a corollary we obtain a result which was proved by Green in the case $R_0 = k$.

**Corollary 2.3** *Let $(R_0, \mathfrak{m})$ be an Artinian local equicharacteristic ring with infinite residue field $k$ and $R = R_0[X_1, \ldots, X_b]$. Let $d \geq 1$ be an integer and $V \subseteq R_d$ a $R_0$-submodule; for all $h \in U_R(d,V)$ we put $\overline{R} = R/(h)$ and $\overline{V} = (V + hR_{d-1})/hR_{d-1}$. Then we have*

$$\lambda_{R_0}(\overline{R}_d/\overline{V}) \leq (\lambda_{R_0}(R_d/V)_d)^-.$$

**Proof.** Let $c \geq 1$ be an integer such that $\lambda_{R_0}(R_d/V) < \binom{d+c-1}{c-1}$.

If $c \leq b$ then $\lambda_{R_0}(R_d/V) < \binom{d+b-1}{b-1}$, so in the Euclidean division of $\lambda_{R_0}(R_0)$ by $\binom{d+b-1}{b-1}$ we get $q = 0$ and $r = \lambda_{R_0}(R_d/V)$. Then by Theorem 2.2 we obtain the result.

Assume $b < c$ and consider $R = R_0[X_1, \ldots, X_b] \subseteq R' = R_0[X_1, \ldots, X_c]$. Define $K = (X_{b+1}, \ldots, X_c)_d$ and $V' = V + K \subseteq R'_d$. Notice that $V \cap K = R_d \cap K = 0$, so we have an isomorphism of $R_0$-modules $R'_d/V' \cong R_d/V$.

Let for $l \in R'_1$, $l = h + f$ with $h \in R_1 \cap k^b$ and $f \in (X_{b+1}, \ldots, X_c)$, then $V' + lR'_{d-1} = V + hR_d + K$ and therefore $l \in U_{R'}(d, V')$ if and only if $h \in U_R(d, V)$. In particular if $h \in U_R(d, V)$ then $h \in U_{R'}(d, V')$. Moreover $\overline{R'_d}/\overline{V'} \cong \overline{R}_d/\overline{V}$, hence

$$\lambda_{R_0}(\overline{R}_d/\overline{V}) = \lambda_{R_0}(\overline{R'}_d/\overline{V'}) \leq (\lambda_{R_0}(\overline{R'}_d/\overline{V'})_d)^- = (\lambda_{R_0}(\overline{R}_d/\overline{V})_d)^-$$

by the former case $c = b$. $\square$

The following proposition is the key tool in the first part of the proof of the main theorem.

**Proposition 2.4** Let $(R_0, \mathfrak{m})$ be an Artinian local equicharacteristic ring, $R = R_0[X_1, \ldots, X_b]$, $d \geq 1$ an integer, $V \subseteq R_d$ a $R_0$-submodule and $\lambda_{R_0}(R_d/V) = \binom{d+b-1}{b-1}q + r$ the Euclidean division. Then
$$\lambda_{R_0}(R_{d+1}/R_1V) \leq \binom{d+1+b-1}{b-1}q + (r_d)_+^+.$$

**Proof.** Notice that without loss of generality we may assume that $R_0$ has infinite residue field. We will proceed by induction on $b$. In the case $b = 1$ we have $\lambda_{R_0}(R_d/V) = q$ and $r = 0$. Since $V \subseteq R_d = R_0\langle X_1^d \rangle$ we have $V = \mathfrak{a} X_1^d$ with $\mathfrak{a} \subseteq R_0$ an ideal. So $R_1 V = \mathfrak{a} X_1^{d+1}$, and then
$$\lambda_{R_0}(R_{d+1}/R_1V) = \lambda_{R_0}(R_0/\mathfrak{a}) = q = q + (0_d)_+^+$$
as we wanted to prove.

In the case $b \geq 2$, pick $h \in U_R(d+1, R_1V)$ and consider the exact sequence
$$R_d/V \xrightarrow{\cdot h} R_{d+1}/R_1V \to R_{d+1}/(R_1V + hR_d) \to 0.$$

Let
$$\lambda_{R_0}(\overline{R}_d/\overline{V}) = \binom{d+b-2}{b-2}\overline{q} + \overline{r}$$
be the Euclidean division. By Theorem 2.2 we have
$$\binom{d+b-2}{b-2}\overline{q} + \overline{r} = \lambda_{R_0}(\overline{R}_d/\overline{V}) \leq \binom{d+b-2}{b-2}q + (r_d)^-$$

and both sides are Euclidean divisions, see the proof of Theorem 2.2. So we have $(\overline{q}, \overline{r}) \leq (q, (r_d)^-)$ in the lexicographical order.

Applying the induction hypothesis to $\overline{R}_d/\overline{V}$ we obtain
$$\lambda_{R_0}(\overline{R}_{d+1}/\overline{R_1V}) = \lambda_{R_0}(\overline{R}_{d+1}/\bar{R}_1 \overline{V}) \leq \binom{d+1+b-2}{b-2}\overline{q} + (\overline{r}_d)_+^+.$$

**Claim:** We have an inequality of Euclidean divisions
$$\binom{d+1+b-2}{b-2}\overline{q} + (\overline{r}_d)_+^+ \leq \binom{d+1+b-2}{b-2}q + ((r_d)^-)_+^+.$$

*Proof:* Recall that $0 \leq \overline{r}, (r_d)^- < \binom{d+b-2}{b-2}$, so we obtain $0 \leq (\overline{r}_d)_+^+, ((r_d)^-)_+^+ < \binom{d+1+b-2}{b-2}$, i.e. the two expressions are Euclidean divisions. On the other hand the inequality $(\overline{q}, \overline{r}) \leq (q, (r_d)^-)$ implies $(\overline{q}, (\overline{r}_d)_+^+) \leq (q, ((r_d)^-)_+^+)$, and this proves the claim.

By the exact sequence we have
$$\begin{aligned}\lambda_{R_0}(R_{d+1}/(R_1V)) &\leq \lambda_{R_0}(R_d/V) + \lambda_{R_0}(\overline{R}_{d+1}/\overline{R_1V}) \\ &\leq \binom{d+b-1}{b-1}q + r + \binom{d+1+b-2}{b-2}q + ((r_d)^-)_+^+ \\ &= \binom{d+1+b-1}{b-1}q + r + ((r_d)^-)_+^+.\end{aligned}$$

From [Rob90], Proposition 4.8, we have $((r_d)^-)_+^+ = (r_d)_+^+ - r$, and so we get the result. $\square$

Let us define the functions we seek to characterize.

**Definition 2.5** *Given an Artinian equicharacteristic local ring $(R_0, \mathfrak{m})$ and a function $H : \mathbb{N} \to \mathbb{N}$, we will say that $H$ is admissible if there exists a standard $R_0$-algebra $S$ such that $H = H_S$. For $b \geq 1$ we will say that $H$ is $b$-admissible if there exists a homogeneous ideal $I \subseteq R_+$, where $R = R_0[X_1, \ldots, X_b]$, such that $H = H_{R/I}$.*

The following result is the main theorem of this section and characterizes under which conditions a function $H$ is $b$-admissible. It is a stronger analogue to Macaulay's theorem. In the case $R_0 = k$ one gets trivially that if $H$ is admissible, then it is $b$-admissible if and only if $b \geq H(1)$, see Corollary 2.11 and Remark 2.12. But in the general case there is a gap for $b$ between $H(1)/\lambda_{R_0}(R_0)$ and $H(1)$ which can only be covered by refining the conditions on $H$.

For the constructive part of the proof we will need to consider an order in $R = R_0[X_1, \ldots, X_b]$ which, in addition to the combinatorics of the monomials, also takes into account the structure of $R_0$. Let us begin by fixing a suitable order in the set of monomials of $R$. Given a multiindex $\lambda = (\lambda_1, \ldots, \lambda_b)$, let $X^\lambda = X_1^{\lambda_1} \ldots X_b^{\lambda_b}$ and $|\lambda| = \lambda_1 + \cdots + \lambda_b$. We have chosen the degree reverse lexicographical order to work with:

**Definition 2.6** *For $X^\lambda$, $X^\mu$ monomials in $R$, we will say that $X^\lambda > X^\mu$ if $|\lambda| > |\mu|$ or $|\lambda| = |\mu|$ and the last nonzero entry of $(\lambda_1 - \mu_1, \ldots, \lambda_b - \mu_b)$ is negative.*

Nevertheless, our constructions would work exactly the same way using the degree lexicographical order instead. Let $\mathcal{J} = \{0 = J_0 \subseteq J_1 \subseteq \ldots \subseteq J_s = R_0\}$, where $s = \lambda_{R_0}(R_0)$, be a composition series in $R_0$ and consider for all $n \geq 1$ the set of $R_0$-submodules of $R_n$

$$\mathcal{M}_n(\mathcal{J}) = \{J_i X^\lambda \mid 1 \leq i \leq r \text{ and } |\lambda| = n\},$$

notice that $\#(\mathcal{M}_n(\mathcal{J})) = s\binom{n+b-1}{b-1}$.

**Definition 2.7** *We define a total ordering, $\mathcal{J}$-reverse lexicographical order, in $\mathcal{M}_n(\mathcal{J})$ by*

$$J_i X^\lambda < J_l X^{\lambda'} \Leftrightarrow i < l \text{ or } i = l \text{ and } X^\lambda < X^{\lambda'}$$

*where the order in the set of monomials in $X_1, \ldots, X_b$ is the degree reverse lexicographical order.*

**Theorem 2.8 (Characterization of Hilbert functions)** *Let $(R_0, \mathfrak{m})$ be an Artinian local equicharacteristic ring, $H : \mathbb{N} \to \mathbb{N}$ a function, $b \geq 1$ and $R = R_0[X_1, \ldots, X_b]$. For all $n \geq 0$ let us consider the Euclidean division*

$$H(n) = \binom{n+b-1}{b-1} q(n) + r(n).$$

*Then the following conditions are equivalent:*

(i) *There exists a homogeneous ideal $I \subseteq R_+$ such that $H = H_{R/I}$,*

(ii) $H(0) = \lambda_{R_0}(R_0)$ *and* $H(n+1) \leq \binom{n+1+b-1}{b-1} q(n) + (r(n)_n)^+_+$ *for all $n \geq 0$.*

**Proof.** Assume that $H$ verifies the condition in (i) and let $s = H(0) = \lambda_{R_0}(R_0)$. For $n = 0$ we have $q(0) = s$ and $r(0) = 0$, and we must show $H(1) \leq bs$. But this is obvious since $H(1) = \lambda_{R_0}(R_1/I_1) \leq \lambda_{R_0}(R_1) = bs$.

For $n \geq 1$ we have $R_1 I_n \subseteq I_{n+1}$ and therefore by Proposition 2.4 we have

$$H(n+1) = \lambda_{R_0}(R_{n+1}/I_{n+1}) \leq \lambda_{R_0}(R_{n+1}/R_1 I_n)$$

$$\leq \binom{n+1+b-1}{b-1} q(n) + (r(n)_n)_+^+.$$

Reciprocally, assume that $H$ verifies the conditions in (ii). Notice that these are equivalent to

$$(0,0) \leq (q(n+1), r(n+1)) \leq (q(n), (r(n)_n)_+^+) \leq (s, 0)$$

in the lexicographical order, for all $n \geq 1$. Let us consider $\mathcal{J} = \{0 = J_0 \subseteq J_1 \subseteq \ldots \subseteq J_s = R_0\}$ a composition series in $R_0$. Fix $n \geq 1$ and let $N = \binom{n+b-1}{b-1}$. Let $X^{\lambda_1} > X^{\lambda_2} > \ldots > X^{\lambda_N}$ be the ordered monomials of degree $n$ in $X_1, \ldots, X_b$ and define the following $R_0$-submodule of $R_n$:

$$I_n = J_{s-q(n)-1} X^{\lambda_1} + \ldots + J_{s-q(n)-1} X^{\lambda_{r(n)}} + J_{s-q(n)} X^{\lambda_{r(n)+1}} + \ldots + J_{s-q(n)} X^{\lambda_N},$$

this makes sense since $(0,0) \leq (q(n), r(n)) \leq (s, 0)$. Since $R_n = R_0 X^{\lambda_1} \oplus \ldots \oplus R_0 X^{\lambda_N}$ we have

$$\frac{R_n}{I_n} \cong \left(\frac{R_0}{J_{s-q(n)-1}}\right)^{r(n)} \oplus \left(\frac{R_0}{J_{s-q(n)}}\right)^{N-r(n)}.$$

On the other hand $\lambda_{R_0}(R_0/J_i) = s - i$, and so we get

$$\lambda_{R_0}(R_n/I_n) = r(n)(s - (s - q(n) - 1)) + (N - r(n))(s - (s - q(n)))$$

$$= r(n)(q(n) + 1) + (N - r(n))q(n)$$

$$= Nq(n) + r(n) = H(n).$$

Hence, to finish the proof it suffices to show that $I = \oplus_{n \geq 1} I_n$ is an ideal of $R$.

Consider for all $n \geq 1$ the ordered set $\mathcal{M}_n(\mathcal{J})$ of $R_0$-submodules of $R_n$ and let us write the elements of $\mathcal{M}_n(\mathcal{J})$ ordered from greater to lesser:

| $J_s X^{\lambda_1}$ | $J_s X^{\lambda_2}$ | $\ldots$ | $J_s X^{\lambda_{r(n)}}$ | $J_s X^{\lambda_{r(n)+1}}$ | $\ldots$ | $J_s X^{\lambda_N}$ |
|---|---|---|---|---|---|---|
| $J_{s-1} X^{\lambda_1}$ | $J_{s-1} X^{\lambda_2}$ | $\ldots$ | $J_{s-1} X^{\lambda_{r(n)}}$ | $J_{s-1} X^{\lambda_{r(n)+1}}$ | $\ldots$ | $J_{s-1} X^{\lambda_N}$ |
| $\vdots$ | $\vdots$ | | $\vdots$ | $\vdots$ | | $\vdots$ |
| $J_{s-q(n)+1} X^{\lambda_1}$ | $J_{s-q(n)+1} X^{\lambda_2}$ | $\ldots$ | $J_{s-q(n)+1} X^{\lambda_{r(n)}}$ | $J_{s-q(n)+1} X^{\lambda_{r(n)+1}}$ | $\ldots$ | $J_{s-q(n)+1} X^{\lambda_N}$ |
| $J_{s-q(n)} X^{\lambda_1}$ | $J_{s-q(n)} X^{\lambda_2}$ | $\ldots$ | $J_{s-q(n)} X^{\lambda_{r(n)}}$ | $J_{s-q(n)} X^{\lambda_{r(n)+1}}$ | $\ldots$ | $J_{s-q(n)} X^{\lambda_N}$ |
| $J_{s-q(n)-1} X^{\lambda_1}$ | $J_{s-q(n)-1} X^{\lambda_2}$ | $\ldots$ | $J_{s-q(n)-1} X^{\lambda_{r(n)}}$ | $J_{s-q(n)-1} X^{\lambda_{r(n)+1}}$ | $\ldots$ | $J_{s-q(n)-1} X^{\lambda_N}$ |
| $\vdots$ | $\vdots$ | | $\vdots$ | $\vdots$ | | $\vdots$ |
| $J_1 X^{\lambda_1}$ | $J_1 X^{\lambda_2}$ | $\ldots$ | $J_1 X^{\lambda_{r(n)}}$ | $J_1 X^{\lambda_{r(n)+1}}$ | $\ldots$ | $J_1 X^{\lambda_N}$, |

read from left to right and from top to bottom.

Notice that $I_n$ can be seen graphically by deleting the first $q(n)$ rows in $\mathcal{M}_n(\mathcal{J})$ and the first $r(n)$ elements in the $(q(n)+1)$-th row, and keeping the remaining elements as generators of $I_n$. In particular, the condition $(q(n), r(n)) \leq (s, 0)$ assures us that we are not trying to delete more rows than we really have.

What we have to prove is $R_1 I_n \subseteq I_{n+1}$. Since $R_1 \cdot$ (i-th row of $\mathcal{M}_n(\mathcal{J})$) $\subseteq$ (i-th row of $\mathcal{M}_{n+1}(\mathcal{J})$) and $q(n+1) \leq q(n)$, in the case $q(n+1) < q(n)$ we have that the generators of $R_1 I_n$ are contained in the generators of $I_{n+1}$. On the other hand, if $q(n) = q(n+1)$ we are deleting the same rows in $\mathcal{M}_n(\mathcal{J})$ and in $\mathcal{M}_{n+1}(\mathcal{J})$, and what we have to prove is then

$$R_1 \cdot (J_{s-q(n)} X^{\lambda_{r(n)+1}} + \cdots + J_{s-q(n)} X^{\lambda_N}) \subseteq J_{s-q(n)} X^{\mu_{r(n+1)+1}} + \cdots + J_{s-q(n)} X^{\mu_M}$$

where $M = \binom{n+1+b-1}{b-1}$. Notice that we can ignore $J_{s-q(n)}$ and it is then enough to show that $X_i \cdot X^{\lambda_{r(n)+1}}, \ldots, X_i \cdot X^{\lambda_N} \in (X^{\mu_{r(n+1)+1}}, \ldots, X^{\mu_M})$ for all $1 \leq i \leq b$. And this, in the same way as in the proof of Macaulay's theorem, is a consequence of the fact $r(n+1) \leq (r(n)_n)^+_+$, see for example [BH93], Proposition 4.2.8. □

The ideal $I \subseteq R_0[X_1, \ldots, X_b]$ verifying $H_{R/I} = H$ has been constructed in the proof of Theorem 2.8 by deleting in each degree $n$ the first $H(n)$ elements in $\mathcal{M}_n(\mathcal{J})$ and taking the remaining ones as generators for $I_n$. Let us give a name to this type of ideals:

**Definition 2.9** *A homogeneous ideal $I \subseteq R = R_0[X_1, \ldots, X_b]_+$ will be called a $\mathcal{J}$-segment ideal if for all $n \geq 1$ $I_n$ is generated as a $R_0$-module by the $s\binom{n+b-1}{b-1} - \lambda_{R_0}(R_n/I_n)$ smallest elements in $\mathcal{M}_n(\mathcal{J})$.*

It turns out that for each $b$-admissible function $H$, there exists a unique $\mathcal{J}$-segment ideal $I \subseteq R_0[X_1, \ldots, X_b]$, such that $H = H_{R/I}$; it will be denoted by $I_{H,\mathcal{J}}$.

Notice that in the case $R_0 = k$ the only composition series is the trivial one, hence the $\mathcal{J}$-reverse lexicographical order and the $\mathcal{J}$-segment ideals coincide with the usual degree reverse lexicographical order and segment ideals. In this case we will write $I_{H,\mathcal{J}} = I_H$.

The results in Section 3 will allow us to effectively apply the characterization theorem; see Section 5 for some examples, in which we check the $b$-admissibility of some functions and compute the ideals $I_{H,\mathcal{J}}$.

The following corollary is the direct generalization of the classical version of Macaulay's theorem:

**Corollary 2.10** *Let $(R_0, \mathfrak{m})$ be an Artinian local equicharacteristic ring and $H : \mathbb{N} \to \mathbb{N}$ a function. Then the following conditions are equivalent:*

*(i) There exists a standard $R_0$-algebra $S$ such that $H = H_S$,*

*(ii) $H(0) = \lambda_{R_0}(R_0)$ and $H(n+1) \leq (H(n)_n)^+_+$ for all $n \geq 1$.*

**Proof.** Let $S = R/I$, where $R = R_0[X_1, \ldots, X_b]$ and $I \subseteq R_+$ is a homogeneous ideal with $H_S = H$. The first part of (ii) is immediate, and for the second part we will distinguish two cases:

Case (a): $H(1) < b$. Let $H(n) = \binom{n+b-1}{b-1} q(n) + r(n)$ be the Euclidean division. Since $bq(1) + r(1) = H(1) < b$ it must be $q(1) = 0$. From Theorem 2.8 we get that $q(n) = 0$ and $r(n) = H(n)$ for all $n \geq 1$, and again by Theorem 2.8 we have (ii).

Case (b): $H(1) \geq b$. Let $c = H(1) + 1 > b$ and $R' = R_0[X_1, \ldots, X_c]$. Consider the natural projection $\pi : R' \to R$ and let $I' = \pi^{-1}(I)$. From $R/I \cong R'/I'$ we get $H = H_{R/I} = H_{R'/I'}$, and we are in case (a).

Reciprocally, assume that $H$ verifies the conditions in (ii) and let $b = H(1)$. From the condition $H(n+1) \leq (H(n)_n)_+^+$ we get for all $n \geq 1$

$$H(n) \leq \binom{n+b-1}{b-1}$$

which is the number of monomials of degree $n$ in $b$ variables. Now let $R = R_0[X_1, \ldots, X_b]$ and consider the monomials of degree $n$, $n \geq 1$, in $X_1, \ldots, X_b$ ordered by the degree reverse lexicographical order: $X^{\lambda_1} > X^{\lambda_2} > \ldots > X^{\lambda_N}$, where $N = \binom{n+b-1}{b-1}$. Consider then the $R_0$-submodule $I_n \subseteq R_n$ defined by

$$I_n = \mathfrak{m}X^{\lambda_1} + \ldots + \mathfrak{m}X^{\lambda_{H(n)}} + R_0 X^{\lambda_{H(n)+1}} + \ldots + R_0 X^{\lambda_N}.$$

Since $R_n = \oplus_{j=1}^N R_0 X^{\lambda_j}$, we have $R_n/I_n = \oplus_{j=1}^{H(n)} R_0/\mathfrak{m}$ and $\lambda_{R_0}(R_n/I_n) = H(n)$. Therefore, the only thing left to complete the proof is to check that $I = \oplus_{n \geq 1} I_n$ is an ideal, i.e. $R_1 I_n \subseteq I_{n+1}$ for all $n$. This is a consequence of the condition $H(n+1) \leq (H(n)_n)_+^+$, see [BH93], Proposition 4.2.8. □

**Corollary 2.11** *Let $H : \mathbb{N} \to \mathbb{N}$ be an admissible function; then $H$ is $b$-admissible for all $b \geq H(1)$.*

**Proof.** By Corollary 2.10 $H$ verifies $H(0) = \lambda_{R_0}(R_0)$ and $H(n+1) \leq (H(n)_n)_+^+$ for all $n \geq 1$. Assume first that $b = H(1)$; in this case, we have constucted in the proof of Corollary 2.10 an ideal $I \subseteq R_+$, where $R = R_0[X_1, \ldots, X_{H(1)}]$, such that $H = H_{R/I}$. Hence $H$ is $H(1)$-admissible.

Assume now $b > H(1)$; then $H(n+1) \leq (H(n)_n)_+^+$ for all $n \geq 1$ implies that $H(n) < \binom{n+b-1}{b-1}$ for all $n \geq 1$. So $q(n) = 0$ and $r(n) = H(n)$ for all $n \geq 1$, and the conditions in Theorem 2.8 (ii) hold trivially in this case. □

**Remark 2.12** *Notice that if $H$ is a $b$-admissible function, it must verify $H(1) = \lambda_{R_0}(R_1/I_1) \leq \lambda_{R_0}(R_1) = bs$, where $s = \lambda_{R_0}(R_0)$. Therefore $H$ will never be $b$-admissible for $b < H(1)/s$.*

*In the case $H(1)/s \leq b < H(1)$, either situation is possible; nevertheless we can assure that if $H$ is $b$-admissible for some $b$, it is also $b'$-admissible for all $b' \geq b$.*

The results above can be applied to the local case: let $(A, \mathfrak{m})$ a local ring and $\mathfrak{a}$ a $\mathfrak{m}$-primary ideal. Consider $R_0 = A/\mathfrak{a}$ and $R/I = gr_\mathfrak{a}(A)$ where $R = R_0[X_1, \ldots, X_b]$ and $b$ is the number of generators of $\mathfrak{a}/\mathfrak{a}^2$ over $A/\mathfrak{a}$. Notice that $H_{R/I} = H_\mathfrak{a}^0$. If $\mathfrak{a} = \mathfrak{m}$ then $R_0 = k$ is obviously equicharacteristic; if $\mathfrak{a} \neq \mathfrak{m}$ we need to assume that $R_0$ is an equicharacteristic ring. For example, if $A$ itself is equicharacteristic this condition is always verified, see [Mat86], Theorem 28.3.

The application of the results of this section to the local case will appear in a forthcoming paper.

# 3  Gotzmann developments of Hilbert polynomials

The main result in this section, Theorem 3.5, is an improved version of Gotzmann's regularity theorem in the Artinian equicharacteristic case, see [Gre89]. This theorem gives an alternative expression of Hilbert polynomials, better suited than the usual one to deal with the combinatorical properties of Hilbert functions. For example, this will allow us to characterize Hilbert polynomials and to encode an entire Hilbert function in a finite amount of data, see Sections 4 and 5. Furthermore, it also provides us with information about the local cohomology of the ring. Let us begin by recalling some facts about local cohomology; see [HIO88], §35, as a reference.

Let $(R_0, \mathfrak{m})$ be an Artinian local ring, $R$ a standard $R_0$-algebra, $M = \bigoplus_{n \geq 0} M_n$ a finitely generated graded $R$-module. We will denote by $H^q_{R_+}(M) = \bigoplus_{n \in \mathbb{Z}} H^q_{R_+}(M)_n$ the $q$-th local cohomology module of $M$ with respect to $R_+$. Since $\operatorname{rad}(R_+) = \mathfrak{m}$, the maximal homogeneous ideal of $R$, we have $H^q_{R_+}(M) = H^q_{\mathfrak{m}}(M)$ for all $q$. It is known that these modules are Artinian and that for all $q$, $n$, $H^q_{R_+}(M)_n$ is a finitely generated $R_0$-module, hence we can define

$$a_q(M) = \min\{n \in \mathbb{N} \mid H^q_{R_+}(M)_n \neq 0\} < +\infty.$$

It is also known that $H^q_{R_+}(M) = 0$ for $q < \operatorname{depth}_{R_+}(M)$ and $q > \dim(M)$. We will adopt the convention that $a_q(M) = -\infty$ for $q < \operatorname{depth}_{R_+}(M)$. The relationship between local cohomology and Hilbert functions is given by the following result:

**Proposition 3.1 (Grothendieck's formula)** *Let $(R_0, \mathfrak{m})$ be an Artinian local ring, $R$ a standard $R_0$-algebra, $M$ a finitely generated graded $R$-module; then*

$$H_M(n) - h_M(n) = \sum_{i \geq 0} (-1)^i \lambda_{R_0}(H^q_{R_+}(M)_n)$$

*for all $n \in \mathbb{Z}$.*

See [Mar93], Lemma 1.3 for a purely algebraic proof.

We will begin by stating three preliminary lemmas.

**Lemma 3.2** *Let $(R_0, \mathfrak{m})$ be an Artinian local equicharacteristic ring with infinite residue field $k$, $R = R_0[X_1, \ldots, X_b]$, $I \subseteq R_+$ a homogeneous ideal such that $\dim(R/I) \geq 1$. The following conditions are equivalent:*

(i) $H^0_{R_+}(R/I) = 0$,

(ii) *There exists $h \in R_1 \cap k^b$ such that $\bar{h}$ is a non-zero divisor in $R/I$.*

**Proof.** If (ii) holds, then in particular $\operatorname{depth}(R/I) \geq 1$ and hence $H^0_{R_+}(R/I) = 0$.

Reciprocally, assume (i) and let $\mathfrak{P}_1, \ldots, \mathfrak{P}_s$ be the associated primes of $R/I$, so $z(R/I) = \mathfrak{P}_1 \cup \ldots \cup \mathfrak{P}_s$. For all $i$ we have $\mathfrak{P}_i = \operatorname{Ann}(\bar{f}_i)$ with $\bar{f}_i \neq 0$, $\bar{f}_i \in R/I$ homogeneous. Since $H^0_{R_+}(R/I) = 0$, we have in particular that $R_+ \bar{f}_i \neq 0$, that is, $R_+ \not\subseteq \mathfrak{P}_i$. Therefore $X_1, \ldots, X_b$ can not simultaneously belong to $\mathfrak{P}_i$, and so $R_1 \cap k^b \not\subseteq (\mathfrak{P}_i)_1 \cap k^b$ for all $1 \leq i \leq s$. In other words, $(\mathfrak{P}_i)_1 \cap k^b$ are proper vector subspaces of $R_1 \cap k^b$. Since $k$ is infinite, we deduce that $((\mathfrak{P}_1)_1 \cap k^b) \cup \ldots \cup ((\mathfrak{P}_s)_1 \cap k^b) \neq R_1 \cap k^b$ and hence we can find a homogeneous element $h \in R_1 \cap k^b$ such that $h \notin \mathfrak{P}_1 \cup \ldots \cup \mathfrak{P}_s$, that is, $\bar{h} \notin z(R/I)$. □

**Lemma 3.3** Let $(R_0, \mathfrak{m})$ be an Artinian local equicharacteristic ring with infinite residue field $k$, $R = R_0[X_1, \ldots, X_b]$ and $I \subseteq R_+$ a homogeneous ideal such that $\text{depth}(R/I) \geq 1$. Then $\bigcap_{n \geq 1} U_R(n, I_n)$ is the set of all non-zero divisors in $R_1 \cap k^b$; in particular it is nonempty.

**Proof.** If $\text{depth}(R/I) \geq 1$, by Lemma 3.2 we may pick $h \in R_1 \cap k^b$ such that $\bar{h}$ is a non-zero divisor in $R/I$. Then $(I_n : h) = I_{n-1}$ for all $n \geq 1$. If $l$ is any element of $R_1 \cap k^b$ we have an exact sequence
$$0 \to (I_n : l) \to R_{n-1} \xrightarrow{\cdot l} (I_n + lR_{n-1})/I_n \to 0.$$
So $\lambda_{R_0}((I_n + lR_{n-1})/I_n)$ is maximal if and only if $\lambda_{R_0}(I_n : l)$ is minimal. Since $I_{n-1} \subseteq (I_n : l)$ and we have equality for all $n \geq 1$ when $l = h$, we get that $l \in \bigcap_{n \geq 1} U_R(n, I_n)$ if and only if $(I_n : l) = I_{n-1}$ for all $n \geq 1$, that is, $\bar{l}$ is not a zero divisor in $R/I$. Hence, $\bigcap_{n \geq 1} U_R(n, I_n)$ is precisely the set of all non-zero divisors in $R_1 \cap k^b$, in particular it contains $h$. □

**Lemma 3.4** Let $(R_0, \mathfrak{m})$ be an Artinian local equicharacteristic ring, $R = R_0[X_1, \ldots, X_b]$, $I \subseteq R_+$ a homogeneous ideal such that $\dim(R/I) \geq 1$. Consider the homogeneous ideal $J \subseteq R$ such that $H^0_{R_+}(R/I) = J/I$ and denote $R'_0 = R_0/J_0$, $R' = R'_0[X_1, \ldots, X_b]$ and $I' = J/J_0 R$. Then we have:

(i) $(R'_0, \mathfrak{m}')$ is an Artinian local equicharacteristic ring.

(ii) $H_{R/I}(n) = H_{R'/I'}(n)$ for all $n \geq a_0(R/I) + 1$. In particular, $R/I$ and $R'/I'$ have the same Hilbert polynomial.

(iii) For all $q \geq 1$, $H^q_{R_+}(R/I) = H^q_{R'_+}(R'/I')$.

(iv) $H^0_{R'_+}(R'/I') = 0$.

**Proof.** (i) The only fact we need to check is that $J_0$ is a proper ideal in $R_0$. If $1 \in J_0$ we would have $R^n_+ \subseteq I$ for some $n$, hence $\dim(R/I) = 0$.

(ii) Since $R' = R'_0[X_1, \ldots, X_b] = R/J_0 R$, there is an isomorphism of graded rings
$$R'/I' \cong \frac{R/J_0 R}{J/J_0 R} \cong R/J,$$
and so $\lambda_{R'_0}(R'_n/I'_n) = \lambda_{R'_0}(R_n/J_n)$ for all $n \geq 0$.

On the other side, by the definition of $a_0$ we have $I_n = J_n$ for all $n \geq a_0 + 1$, that is, $H_{R/I}(n) = \lambda_{R_0}(R_n/I_n) = \lambda_{R_0}(R_n/J_n)$ for all $n \geq a_0 + 1$.

Since $J_0 \subseteq \text{Ann}_{R_0}(R_n/J_n)$, the submodule lattices of $R_n/J_n$ considered as $R_0$ or $R'_0$-module are the same, and so $\lambda_{R_0}(R_n/J_n) = \lambda_{R'_0}(R_n/J_n) = H_{R'/I'}(n)$ for all $n \geq 0$.

(iii) Consider the exact sequence of graded $R$-modules
$$0 \to J/I \to R/I \to R'/I' \to 0.$$

We know that $J/I$ is 0 in big enough degrees, in particular $\dim(J/I) = 0$ and so $H^0_{R_+}(J/I) = J/I$ and $H^q_{R_+}(J/I) = 0$ for $q \geq 1$. Applying it to the local cohomology long exact sequence we obtain for all $q \geq 1$ $H^q_{R_+}(R/I) \cong H^q_{R_+}(R'/I')$, and since the structure of $R$-module of $R'/I'$ is induced by the morphism $R \to R' = R/J_0 R$ we get $H^q_{R_+}(R'/I') \cong H^q_{R'_+}(R'/I')$, see [HIO88], Corollary 35.20.

(iv) We also have by the first part of the local cohomology long exact sequence

$$0 \to J/I \to H^0_{R_+}(R/I) \to H^0_{R_+}(R'/I') \to 0,$$

and the first morphism is an isomorphism, so $H^0_{R_+}(R'/I') = 0$. □

We are ready now to prove the main theorem in this section:

**Theorem 3.5** *Let $(R_0, \mathfrak{m})$ be an Artinian local equicharacteristic ring, $R = R_0[X_1, \ldots, X_b]$ and $I \subseteq R_+$ a homogeneous ideal such that $\dim(R/I) \geq 1$. Then:*

(i) *There exist integers $b - 1 > c'_1 \geq \ldots \geq c'_p \geq 0$, $p \geq 0$, and $0 \leq q \leq \lambda_{R_0}(R_0)$ such that for all $n \gg 0$ we have*

$$h_{R/I}(n) = q\binom{n+b-1}{b-1} + \binom{n+c'_1}{c'_1} + \binom{n+c'_2-1}{c'_2} + \cdots + \binom{n+c'_p-(p-1)}{c'_p};$$

*this equality is an Euclidean division and it holds for all $n \geq \max\{p-1, 0\}$.*

(ii) *Let $p_i = \#\{j \mid c'_j \geq i - 1\}$; then for $i \geq 1$*

$$H^i_{R_+}(R/I)_n = 0 \quad \begin{cases} \text{for all} \quad n \geq p_i - i + 1 & \text{if } q > 0 \\ \text{for all} \quad n \geq p_i - i & \text{if } q = 0. \end{cases}$$

(iii) *The regularity index of $R/I$ verifies*

$$i(R/I) \leq \begin{cases} \max\{a_0(R/I) + 1, p\} & \text{if } q > 0 \\ \max\{a_0(R/I) + 1, p - 1\} & \text{if } q = 0. \end{cases}$$

**Proof.** We may assume without loss of generality that $R_0$ has infinite residue field. Let $r = \lambda_{R_0}(R_0)$. By Lemma 3.4 we may assume that $H^0_{R_+}(R/I) = 0$ (notice that $\lambda_{R'_0}(R'_0) = \lambda_{R_0}(R'_0) \leq \lambda_{R_0}(R_0)$), and in (iii) we have to prove then

$$i(R/I) \leq \begin{cases} p & \text{if } q > 0 \\ p - 1 & \text{if } q = 0. \end{cases}$$

We will proceed by induction on $b$.

In the case $b = 1$ we have $R = R_0[X_1]$. Since $I \subseteq R_+$ is homogeneous we have $I = (\alpha_1 X_1^{m_1}, \ldots, \alpha_l X_1^{m_l})$ with $\alpha_i \in R_+$, $0 < m_1 \leq \ldots \leq m_l$. Then $\alpha_1, \ldots, \alpha_l \in H^0_{R_+}(R/I) = 0$, so $I = 0$.

We must prove (i), i.e.

$$h_R(n) = q\binom{n+b-1}{b-1}$$

for all $n \gg 0$. Since $H_R(n) = \lambda_{R_0}(R_n) = \lambda_{R_0}(R_0 X_1^n) = \lambda_{R_0}(R_0) = r$ for all $n$, (i) holds true taking $q = r > 0$ and $p = 0$. Besides $i(R) = 0 = p$, so it only remains to show that $H^1_{R_+}(R/I)_n = 0$ for all $n \geq 0 = p$.

By Grothendieck's formula 3.1 we have for all $n \geq 0$

$$0 = H_R(n) - h_R(n) = \sum_{i \geq 0}(-1)^i \lambda_{R_0}(H^i_{R_+}(R)_n) = -\lambda_{R_0}(H^1_{R_+}(R)_n)$$

so we obtain (iii) in the case $b = 1$.

In the case $b \geq 2$, by Lemma 3.2 we can choose $h \in R_1 \cap k^b$ such that $\bar{h} \notin z(R/I)$. Let $S = R/(h)$ and $J = (I + (h))/(h)$; we have an exact sequence of graded $R$-modules

$$(*) \quad 0 \to R/I(-1) \xrightarrow{\cdot h} R/I \to S/J \to 0.$$

If $\dim(S/J) = 0$ then $\dim(R/I) = 1$ and $h_{R/I}$ has degree 0, say $h_{R/I} = t$. Since $b - 1 > 1$, (i) holds taking $q = 0$, $p = t$ and $c'_1 = \cdots = c'_p = 0$. We need to show $H^1_{R_+}(R/I)_n = 0$ for all $n \geq p - 1$ and $i(R/I) \leq p - 1$; the last one is a well-known fact since $R/I$ is Cohen-Macaulay. Then, again from Grothendieck's formula we obtain

$$H_{R/I}(n) - h_{R/I}(n) = -\lambda_{R_0}(H^1_{R_+}(R/I)_n)$$

for all $n \geq 0$, and hence $H^1_{R_+}(R/I)_n = 0$ for all $n \geq p - 1$.

Assume now $\dim(S/J) \geq 1$. Consider $S'$ and $J'$ obtained from $S$ and $J$ as in Lemma 3.4; we have $\operatorname{rank}_{R'_0}(S') = \operatorname{rank}_{R_0}(S) = b - 1$ and $\lambda_{S'_0}(S'_0) = \lambda_{R_0}(S'_0) \leq \lambda_{R_0}(S_0) = \lambda_{R_0}(R_0) = r$. So the induction hypothesis applies to $S'/J'$ because of 3.4 (iv), and so by 3.4 (ii)

$$h_{S/J}(n) = h_{S'/J'}(n) = q\binom{n+b-2}{b-2} + \binom{n+b'_1}{b'_1} + \binom{n+b'_2-1}{b'_2} + \cdots + \binom{n+b'_v-(v-1)}{b'_v}$$

for all $n \gg 0$, with $0 \leq q \leq \lambda_{S'_0}(S'_0) \leq \lambda_{R_0}(R_0)$ and $b - 2 > b'_1 \geq \ldots \geq b'_v \geq 0$. Moreover, by 3.4 (iii)

$$H^i_{S_+}(S/J)_n \cong H^i_{S'_+}(S'/J')_n = 0$$

for all $i \geq 1$, $n \geq v_i - i$, where $v_i = \#\{j \mid b'_j \geq i - 1\}$.

To prove (i), fix $n_0 \geq i(R/I), i(S'/J'), a_0(S/J) + 1, v$. Then, since by $(*)$ we have for all $n \geq 0$ that $H_{R/I}(n) - H_{R/I}(n-1) = H_{S/J}(n)$, by Lemma 3.4 (ii) we obtain for all $n \geq n_0$

$$\begin{aligned}
h_{R/I}(n) &= \sum_{j=n_0}^n h_{S/J}(j) + H_{R/I}(n_0 - 1) \\
&= \sum_{j=n_0}^n q\binom{j+b-2}{b-2} + \sum_{j=n_0}^n \sum_{i=1}^v \binom{j+b'_i-(i-1)}{b'_i} + H_{R/I}(n_0 - 1) \\
&= q\binom{n+b-1}{b-1} - q\binom{n_0+b-2}{b-1} + \sum_{i=1}^v \sum_{j=i-1}^n \binom{j+b'_i-(i-1)}{b'_i} - \sum_{i=1}^v \sum_{j=i-1}^{n_0-1} \binom{j+b'_i-(i-1)}{b'_i} \\
&\quad + H_{R/I}(n_0 - 1) \\
&= q\binom{n+b-1}{b-1} - q\binom{n_0+b-2}{b-1} + \sum_{i=1}^v \binom{j+b'_i+1-(i-1)}{b'_i+1} - \sum_{i=1}^v \binom{n_0+b'_i-(i-1)}{b'_i+1} + H_{R/I}(n_0 - 1) \\
&= q\binom{n+b-1}{b-1} + \binom{n+c'_1}{c'_1} + \binom{n+c'_2-1}{c'_2} + \cdots + \binom{n+c'_v-(v-1)}{c'_v} + \rho
\end{aligned}$$

where $c'_i = b'_i + 1$ and $\rho = H_{R/I}(n_0 - 1) - \sum_{i=1}^v \binom{n_0+b'_i-(i-1)}{b'_i+1} - q\binom{n_0+b-2}{b-1}$ is an integer independent of $n$.

If $\rho \geq 0$, taking $p = v + \rho$ and $c'_{v+1} = \cdots = c'_p = 0$ we will have (i). Assume then $\rho < 0$; then for $n \gg 0$ we would have

$$H_{R/I}(n) < q\binom{n+b-1}{b-1} + \binom{n+c'_1}{c'_1} + \binom{n+c'_2-1}{c'_2} + \cdots + \binom{n+c'_v-(v-1)}{c'_v}$$

$$= q\binom{n+b-1}{b-1} + \binom{n+c'_1}{n} + \binom{n+c'_2-1}{n-1} + \cdots + \binom{n+c'_v-(v-1)}{n-(v-1)}.$$

Notice that for $n \geq v$ this is an Euclidean division since $c'_1 < b - 1$ and we have the $n$-binomial expansion of the remainder. By Lemma 3.3, since $h \notin z(R/I)$, we can apply Theorem 2.2 to get for all $n \gg 0$

$$H_{S/J}(n) < q\binom{n+b-2}{b-2} + \binom{n+c'_1-1}{n} + \binom{n+c'_2-1-1}{n-1} + \cdots + \binom{n+c'_v-1-(v-1)}{n-(v-1)}$$

$$= q\binom{n+b-2}{b-2} + \binom{n+b'_1}{b'_1} + \binom{n+b'_2-1}{b'_2} + \cdots + \binom{n+b'_v-(v-1)}{b'_v} = h_{S/J}(n)$$

the strict inequality being consequence of the fact $c'_v = b'_v + 1 \geq 1$ together with Lemma 2.1 (ii).

Now to prove (ii) let us observe that by the definition of the $c'_i$ we have for all $i \geq 2$ that $p_i = v_{i-1}$. Let $q > 0$ (resp. $q = 0$).

By the local cohomology long exact sequence associated to $(*)$ we have for all $n$ and for all $i \geq 2$

$$\cdots \to H^{i-1}_{R_+}(S/J)_n \to H^i_{R_+}(R/I)_{n-1} \to H^i_{R_+}(R/I)_n \to H^i_{R_+}(S/J)_n \to \cdots$$

Since $H^i_{R_+}(S/J)_n = H^i_{S_+}(S/J)_n = 0$ for all $n \geq v_i - i + 1$ (resp. $n \geq v_i - i$) and $v_{i-1} \geq v_i$, we get

$$H^i_{R_+}(R/I)_{n-1} \cong H^i_{R_+}(R/I)_n$$

for all $n \geq v_{i-1} - (i-1) + 1$ (resp. $n \geq v_{i-1} - (i-1)$). Hence $H^i_{R_+}(R/I)_n = 0$ for all $n \geq v_{i-1} - i + 1 = p_i - i + 1$ (resp. $n \geq v_{i-1} - i = p_i - i$).

The only thing left to complete the proof is to show that $H^1_{R_+}(R/I)_n = 0$ for all $n \geq p$ (resp. $n \geq p - 1$), and $H_{R/I}(n) = h_{R/I}(n)$ for all $n \geq p$ (resp. $n \geq p - 1$). Notice that $p = p_1 \geq p_2 \geq \ldots$, so we have just seen that for $n \geq p$ (resp. $n \geq p - 1$) $H^i_{R_+}(R/I)_n = 0$ for all $i \geq 2$. Therefore by Grothendieck's formula we get

$$H_{R/I}(n) - h_{R/I}(n) = -\lambda_{R_0}(H^1_{R_+}(R/I)_n)$$

for all $n \geq p$ (resp. $n \geq p - 1$), that is, for all $n \geq p$ (resp. $p - 1$) we have $H_{R/I}(n) \leq h_{R/I}(n)$, with equality if and only if $H^1_{R_+}(R/I)_n = 0$.

Assume $H^1_{R_+}(R/I)_n \neq 0$ for some $n \geq p$, then we would have

$$H_{R/I}(n) < h_{R/I}(n) = q\binom{n+b-1}{b-1} + \binom{n+c'_1}{c'_1} + \binom{n+c'_2-1}{c'_2} + \cdots + \binom{n+c'_p-(p-1)}{c'_p}$$

$$= q\binom{n+b-1}{b-1} + \binom{n+c'_1}{n} + \binom{n+c'_2-1}{n-1} + \cdots + \binom{n+c'_p-(p-1)}{n-(p-1)}.$$

Since this is an Euclidean division, repeatedly applying Theorem 2.8 we get for all $i \geq n$

$$H_{R/I}(i) < q\binom{i+b-1}{b-1} + \binom{i+c'_1}{i} + \binom{i+c'_2-1}{i-1} + \cdots + \binom{i+c'_p-(p-1)}{i-(p-1)}$$

$$= q\binom{i+b-1}{b-1} + \binom{i+c'_1}{c'_1} + \binom{i+c'_2-1}{c'_2} + \cdots + \binom{i+c'_p-(p-1)}{c'_p} = h_{R/I}(i),$$

contradicting the definition of $h_{R/I}$. Thus we get the result in the case $q > 0$, and in the case $q = 0$ it only remains to show that $H^1_{R_+}(R/I)_{p-1} = 0$. Again, if $H^1_{R_+}(R/I)_{p-1} \neq 0$ we would have

$$H_{R/I}(p-1) < \binom{p-1+c'_1}{c'_1} + \binom{p-1+c'_2-1}{c'_2} + \cdots + \binom{p-1+c'_p-(p-1)}{c'_p}$$

$$= \binom{p-1+c'_1}{p-1} + \binom{p-1+c'_2-1}{p-2} + \cdots + \binom{p-1+c'_p-(p-1)}{0}$$

that is,

$$H_{R/I}(p-1) \leq \binom{p-1+c'_1}{p-1} + \binom{p-1+c'_2-1}{p-2} + \cdots + \binom{p-1+c'_{p-1}-(p-2)}{1}$$

and this is a $(p-1)$-binomial expansion. By Theorem 2.8 we get for all $n \geq p-1$

$$H_{R/I}(n) \leq \binom{n+c'_1}{n} + \binom{n+c'_2-1}{n-1} + \cdots + \binom{n+c'_{p-1}-(p-2)}{n-(p-2)}$$

$$= \binom{n+c'_1}{c'_1} + \binom{n+c'_2-1}{c'_2} + \cdots + \binom{n+c'_{p-1}-(p-2)}{c'_{p-1}} < h_{R/I}(n)$$

and this contradicts the definition of $h_{R/I}$. □

**Remark 3.6** *Notice that in the case $R/I = gr_{\mathfrak{a}}(A)$ with $depth(A) \geq 1$, by [Hoa93a], Theorem 5.2, we have $a_0(R/I) < a_1(R/I)$. Hence we can assure that the maximum in Theorem 3.5 (iii) is $p-1$ if $q=0$ (resp. $p$ if $q>0$).*

As a corollary we obtain a result which was proved by Gotzmann and Green in the case $R_0 = k$.

**Corollary 3.7** *Let $(R_0, \mathfrak{m})$ be an Artinian local equicharacteristic ring, $R = R_0[X_1, \ldots, X_b]$ and $I \subseteq R_+$ a homogeneous ideal such that $dim(R/I) \geq 1$. Then:*

*(i) There exist integers $c_1 \geq \ldots \geq c_s \geq 0$ such that*

$$h_{R/I}(n) = \binom{n+c_1}{c_1} + \binom{n+c_2-1}{c_2} + \cdots + \binom{n+c_s-(s-1)}{c_s}$$

*for all $n \geq s-1$.*

*(ii) Let $s_i = \#\{j \mid c_j \geq i-1\}$; then for $i \geq 1$*

$$H^i_{R_+}(R/I)_n = 0 \quad \forall\, n \geq s_i - i.$$

*(iii) $i(R/I) \leq max\{a_0(R/I) + 1, s-1\}$.*

**Proof.** Let $d = \dim(R/I) \leq b$; since the degree of the polynomial $h_{R/I}$ is $d-1$, if $d < b$ we get $q = 0$ and the statement follows taking $s = p$ and $c_i = c'_i$ for all $1 \leq i \leq s$.

So, we may assume $\dim(R/I) = b$. Notice that in this case we have $height(I) = 0$. Let $R' = R_0[X_1, \ldots, X_{b+1}]$, $\pi : R' \to R$ the projection and $I' = \pi^{-1}(I)$. Since $R/I \cong R'/I'$ we have $H^i_{R_+}(R/I) \cong H^i_{R'_+}(R'/I')$ for all $i$ and $H_{R/I} = H_{R'/I'}$; in particular $h_{R/I} = h_{R'/I'}$. But $R'/I'$ verifies $\dim(R'/I') = b < \dim(R') = b+1$, so we have reduced the problem to the former case. □

**Remark 3.8** In Theorem 3.5 we have in fact that $q \neq 0$ if and only if $\text{height}(I) = 0$ (that is, $\dim(R/I) = \dim(R) = b$). So the case where Theorem 3.5 and Corollary 3.7 are different is the case $R_0 \neq k$ and $I \subseteq R_+$ an ideal contained in $\mathfrak{m}R$, which is the only minimal prime in $R$. In Section 4 we will derive combinatorical consequences from the relationship between these two results and we will check that in fact we always have $p_i \leq s_i$, see Proposition 4.13.

**Proposition 3.9** Let $R_0$ be an Artinian equicharacteristic local ring, $b \geq 1$ an integer, $R = R_0[X_1, \ldots, X_b]$ and $I \subseteq R_+$ a homogeneous ideal, and consider the development

$$h_{R/I}(X) = q\binom{X+b-1}{b-1} + \binom{X+c'_1}{c'_1} + \binom{X+c'_2-1}{c'_2} + \cdots + \binom{X+c'_p-(p-1)}{c'_p}$$

as in Theorem 3.5. Let $m = \max\{i(R/I), p\}$ and for $n \geq 0$ let us consider the Euclidean division $H(n) = \binom{n+b-1}{b-1}q(n) + r(n)$ as in Theorem 2.8. Then we have:

(i) For all $n \geq m$ it holds $q(n) = q$ and $r(n+1) = (r(n)_n)^+_+$.

(ii) $I$ is generated in degrees at most $m$ and in particular if $H^0_{R_+}(R/I) = 0$ then $I$ is generated in degrees at most $p$.

**Proof.** For all $n \geq m$ we have

$$H_{R/I}(n) = h_{R/I}(n) = q\binom{n+b-1}{b-1} + \binom{n+c'_1}{c'_1} + \binom{n+c'_2-1}{c'_2} + \cdots + \binom{n+c'_p-(p-1)}{c'_p}.$$

Since by Theorem 3.5 this is an Euclidean division we must have $q(n) = q$ and

$$r(n) = \binom{n+c'_1}{c'_1} + \binom{n+c'_2-1}{c'_2} + \cdots + \binom{n+c'_p-(p-1)}{c'_p},$$

so (i) is clear. To prove (ii), notice that since $R_1 I_n \subseteq I_{n+1}$ we have $H_{R/I}(n+1) = \lambda_{R_0}(R_{n+1}/I_{n+1}) \leq \lambda_{R_0}(R_{n+1}/R_1 I_n)$. By Proposition 2.4, $\lambda_{R_0}(R_{n+1}/R_1 I_n) \leq \binom{n+b-1}{b-1}q(n) + (r(n)_n)^+_+$. If $n \geq m$, then by (i) $\binom{n+b-1}{b-1}q(n) + (r(n)_n)^+_+ = H_{R/I}(n+1)$, therefore all the inequalities must be equalities and in particular $I_{n+1} = R_1 I_n$ for all $n \geq m$.

Notice also that if $H^0_{R_+}(R/I) = 0$ then $a_0(R/I) = -\infty$ and hence by Theorem 3.5 (iii) we have $i(R/I) \leq p$. $\square$

**Definition 3.10** Let $P \in \mathbb{Q}[X; \mathbb{N}]$; we will say that $P$ admits a Gotzmann development if either $P = 0$ or there exist integers $c_1 \geq c_2 \geq \ldots c_s \geq 0$ such that

$$P(n) = \binom{n+c_1}{c_1} + \binom{n+c_2-1}{c_2} + \cdots + \binom{n+c_s-(s-1)}{c_s}$$

for all $n \gg 0$.

In this case we will call the expression above the Gotzmann development of $P$. Notice that $c_1, \ldots, c_s$ are uniquely determined by $P$; they will be called the Gotzmann coefficients of $P$. We define also $s_q = \#\{i \mid c_i \geq q-1\}$ for all $q \geq 1$.

The first fact to notice is that not all polynomials in $\mathbb{Q}[X; \mathbb{N}]$ admit a Gotzmann development. An example is $P(X) = 2X$. Indeed, if $2n = \binom{n+c_1}{c_1} + \binom{n+c_2-1}{c_2} + \cdots + \binom{n+c_s-(s-1)}{c_s}$ for all $n \gg 0$, it should be $c_1 = c_2 = 1$ and $c_3 = \cdots = c_s = 0$. But this is impossible, because then we would have
$$2n = \binom{n+1}{1} + \binom{n+1-1}{1} + (s-2) = 2n + 1 + s - 2 = 2n + s - 1$$
which is absurd since $s = s_1 \geq s_2 = 2$.

Next we give an equation system to compute the normalized Hilbert-Samuel coefficients from the Gotzmann coefficients and reciprocally:

**Proposition 3.11** *Let $P \in \mathbb{Q}[X; \mathbb{N}]$ a polynomial of degree $d-1$, $d \geq 1$, admitting a Gotzmann development. Denote by $c_1, \ldots, c_s$ the Gotzmann coefficients of $P$ and by $e_0, \ldots, e_{d-1}$ the normalized Hilbert-Samuel coefficients of $P$. Then we have:*

$$\begin{aligned}
e_0 &= s_d \\
e_1 &= \binom{s_d+1}{2} - s_{d-1} \\
e_2 &= \binom{s_d+1}{3} - \binom{s_{d-1}+1}{2} + s_{d-2} \\
&\cdots \\
e_{d-1} &= \binom{s_d+1}{d} - \binom{s_{d-1}+1}{d-1} + \cdots + (-1)^{d-2}\binom{s_2+1}{2} + (-1)^{d-1}s_1
\end{aligned}$$

**Proof.** By induction on $d$. If $d = 1$, $P(n) = e_0$ for all $n$, and its Gotzmann coefficients are $c_1 = \cdots = c_{e_0} = 0$, so $s = s_1 = e_0$. In the case $d \geq 1$ we have
$$P(X) = e_0 \binom{X+d-1}{d-1} - e_1 \binom{X+d-2}{d-2} + \cdots + (-1)^{d-2} e_{d-2} \binom{X+1}{1} + (-1)^{d-1} e_{d-1}$$
for all $X$ and
$$P(X) = \binom{X+c_1}{c_1} + \binom{X+c_2-1}{c_2} + \cdots + \binom{X+c_s-(s-1)}{c_s}$$
for all $X \geq s - 1$.

Given any function $f : \mathbb{Z} \to \mathbb{Z}$ we define $\Delta f(n) = f(n) - f(n-1)$ for all $n \in \mathbb{Z}$. Notice that again $\Delta f : \mathbb{Z} \to \mathbb{Z}$, $\Delta$ is a $\mathbb{Z}$-linear operator and $\Delta \binom{X+i}{i} = \binom{X+i-1}{i-1}$. We have then
$$\Delta P(X) = e_0 \binom{X+d-2}{d-2} - e_1 \binom{X+d-3}{d-3} + \cdots + (-1)^{d-2} e_{d-2}$$
for all $X$, and
$$\Delta P(X) = \binom{X+b_1}{b_1} + \binom{X+b_2-1}{b_2} + \cdots + \binom{X+b_t-(t-1)}{b_t}$$
for all $X \geq s-1$, where $t = s_2$ and $b_i = c_i - 1$. If $t_q = \#\{i \mid b_i \geq q-1\}$, by induction hypothesis

we get

$$e_0 = t_{d-1}$$

$$e_1 = \binom{t_{d-1}+1}{2} - t_{d-2}$$

$$e_2 = \binom{t_{d-1}+1}{3} - \binom{t_{d-2}+1}{2} + t_{d-3}$$

$$\ldots$$

$$e_{d-2} = \binom{t_{d-1}+1}{d-1} - \binom{t_{d-2}+1}{d-2} + \cdots + (-1)^{d-3}\binom{t_2+1}{2} + (-1)^{d-2}t_1.$$

Since $t_q = s_{q+1}$ for all $q \geq 1$ it only remains to compute $e_{d-1}$.

Let us give before an expression of $P$ which involves $s_q$:

$$P(X) = \sum_{i=1}^{s} \binom{X+c_i-(i-1)}{c_i}$$

$$= \sum_{j=1}^{d} \sum_{i=s_{j+1}+1}^{s_j} \binom{X+j-1-(i-1)}{j-1}$$

$$= \sum_{j=1}^{d} \left( \binom{X+j-s_{j+1}}{j} - \binom{X+j-1-(s_j-1)}{j} \right)$$

$$= \binom{X+d}{d} - \sum_{j=2}^{d} \left( \binom{X+j-1-(s_j-1)}{j} - \binom{X+j-1-s_j}{j-1} \right) - \binom{X-s_1+1}{1}$$

$$= \binom{X+d}{d} - \sum_{j=2}^{d} \binom{X+j-s_j-1}{j} - \binom{X-s_1+1}{1},$$

note that $s_{d+1} = 0$.

Comparing the two expressions we have for $P$ and evaluating at $X = -1$ we get

$$(-1)^{d-1}e_{d-1} = -\sum_{j=2}^{d} \binom{j-s_j-2}{j} - \binom{-s_1}{1}$$

$$= -(-s_1) - \sum_{j=2}^{d} \frac{(j-s_j-2)(j-s_j-3)\ldots(-s_j)(-s_j-1)}{j!}$$

$$= s_1 - \sum_{j=2}^{d} (-1)^j \frac{(s_j-j+2)(s_j-j+3)\ldots s_j(s_j+1)}{j!}$$

$$= s_1 + \sum_{j=2}^{d} (-1)^{j+1} \binom{s_j+1}{j},$$

so we obtain

$$e_{d-1} = (-1)^{d-1}s_1 + \sum_{j=2}^{d}(-1)^{d-j}\binom{s_j+1}{j}$$

as we wanted to prove. □

**Remark 3.12** Let us recall a result due to Mumford: For all $b \geq 1$ there exists a polynomial $F_b(t_0, \ldots, t_{b-1})$ such that for all coherent sheaves $\mathcal{I}$ on $\mathbb{P}^{b-1}$, if $a_0, \ldots, a_{b-1}$ are defined by

$$\chi(\mathcal{I}(n)) = a_0\binom{n}{0} + a_1\binom{n}{1} + \cdots + a_{b-1}\binom{n}{b-1}$$

then $\mathcal{I}$ is $F_b(a_0, \ldots, a_{b-1})$-regular. See Lecture 14 in [Mum66].

Now we can make Mumford's result effective: Proposition 3.11 allows us to compute $F_b(a_0, \ldots, a_{b-1})$ for all ideal sheaves $\mathcal{I}$. Indeed, first we compute $e_0, \ldots, e_{d-1}$ via the polynomial identity

$$\binom{n+b-1}{b-1} - \left(a_0\binom{n}{0} + \cdots + a_{b-1}\binom{n}{b-1}\right) = e_0\binom{n+d-1}{d-1} - \cdots + (-1)^{d-1}e_{d-1},$$

*notice that $e_0, \ldots e_{d-1}$ are polynomical on $a_0, \ldots a_{b-1}$. Then we compute $s = s_1$, which is polynomical on $e_0, \ldots, e_{d-1}$. According to Gotzmann's regularity theorem, s is a bound for the regularity of $\mathcal{I}$.*

**Remark 3.13** *The bounds obtained in Corollary 3.7*

$$a_i(R/I) \leq s_i - i - 1$$

*are polynomical functions on $e_0, \ldots, e_d$. For instance, we obtain $a_d(R/I) \leq e_0 - d - 1$ (which is obtained for any Artinian ring $R_0$ in [Hoa93a], Lemma 4.2) and $a_{d-1}(R/I) \leq \binom{e_0+1}{2} - e_1 - d$.*

**Corollary 3.14** *For all $i \geq 1$ we obtain bounds*

$$(-1)^i e_i \geq f_i(e_0, \ldots, e_{i-1}).$$

**Proof.** We have
$$\begin{aligned}
(-1)^i e_i &= (-1)^i \binom{s_d+1}{i+1} + (-1)^{i-1}\binom{s_{d-1}+1}{i} + \cdots + \binom{s_{d-i+2}+1}{3} - \binom{s_{d-i+1}+1}{2} + s_{d-i} \\
&= (-1)^i \binom{s_d+1}{i+1} + (-1)^{i-1}\binom{s_{d-1}+1}{i} + \cdots + \binom{s_{d-i+2}+1}{3} - \binom{s_{d-i+1}}{2} - s_{d-i+1} + s_{d-i} \\
&\geq (-1)^i \binom{s_d+1}{i+1} + (-1)^{i-1}\binom{s_{d-1}+1}{i} + \cdots + \binom{s_{d-i+2}+1}{3} - \binom{s_{d-i+1}}{2} \\
&= f_i(e_0, \ldots, e_{i-1}),
\end{aligned}$$

using that $s_d, s_{d-1}, \ldots, s_{d-i+1}$ depend only on $e_0, \ldots, e_{i-1}$. □

**Remark 3.15** *Notice that for $i = 1$ we get*

$$e_1 \leq \binom{e_0}{2}.$$

*In the local case, $e_i = e_i(\mathfrak{a})$ with $\mathfrak{a}$ a $\mathfrak{m}$-primary ideal such $A/\mathfrak{a}$ is equicharacteristic, the bound for the maximal ideal was first given by G. Valla using Gotzmann's result in an unpublished proof. We learnt it from a private communication of J. Elias, which arose our interest in the subject. L. T. Hoa also proved this result in [Hoa93b], Lemma 5.1. The bound in the case $\mathfrak{a} \neq \mathfrak{m}$ answers a question he posed in that paper in the case $A/\mathfrak{a}$ equicharacteristic.*

*Notice also the the bounds for the higher $e_i$ become much more complicated, for example we obtain*

$$e_2 \geq \binom{e_0+1}{3} - \binom{\binom{e_0+1}{2} - e_1}{2}$$

$$e_3 \leq \binom{e_0+1}{4} - \binom{\binom{e_0+1}{2} - e_1 + 1}{3} + \binom{\binom{\binom{e_0+1}{2} - e_1 + 1}{2} - \binom{e_0+1}{3} + e_2}{2}.$$

# 4 Characterization of Hilbert polynomials

This section is devoted to studying under which conditions a polynomial $P \in \mathbb{Q}[X; \mathbb{N}]$ can be the Hilbert polynomial of a standard algebra. We will see that the good strategy to solve the problem is to use Gotzmann developments; it is shown in Theorem 4.4 that $P$ is a Hilbert polynomial if and only if it admits a Gotzmann development. We also characterize the minimal number of variables for which $P$ is admissible. Let us remark that all the characterizations given are effective; we describe the algorithms in Section 5.

**Definition 4.1** *Given an Artinian equicharacteristic local ring $(R_0, \mathfrak{m})$ and $P \in \mathbb{Q}[X; \mathbb{N}]$, we will say that $P$ is admissible if it is the Hilbert polynomial of a standard $R_0$-algebra. If $b \geq 1$ is an integer, we will say that $P$ is $b$-admissible if there exists a homogeneous ideal $I \subseteq R_+$, where $R = R_0[X_1, \ldots, X_b]$, such that $P = h_{R/I}$.*

Notice that if $P$ is $b$-admissible then $P$ is $b'$-admissible for all $b' \geq b$, see Remark 2.12.

In order to decide whether $P$ is an admissible polynomial it will suffice to study under which conditions $P$ can be interpolated by an admissible function. Let us define a special admissible function that will do:

**Definition 4.2** *Let $c_1 \geq c_2 \geq \ldots c_s \geq 0$ be integers. We define the Gotzmann function $G[c_1, \ldots, c_s]$ associated to $c_1, \ldots, c_s$, by*

$$G[c_1, \ldots, c_s](n) = \begin{cases} \binom{n+c_1}{c_1} + \binom{n+c_2-1}{c_2} + \cdots + \binom{n+c_{n+1}-n}{c_{n+1}} & \text{if } n \leq s-1 \\ \binom{n+c_1}{c_1} + \binom{n+c_2-1}{c_2} + \cdots + \binom{n+c_s-(s-1)}{c_s} & \text{if } n \geq s-1. \end{cases}$$

**Lemma 4.3** *For all $c_1 \geq c_2 \geq \ldots \geq c_s \geq 0$ integers, $G[c_1, \ldots, c_s]$ is an admissible function such that $G[c_1, \ldots c_s](1) = c_1 + 2$.*

**Proof.** Let $H = G[c_1, \ldots, c_s]$. Since for $n \geq s$ the expression of $H(n)$ coincides with its $n$-binomial expansion, we get that $H(n+1) = (H(n)_n)_+^+$ for all $n \geq s$. For $n \leq s-1$ we have

$$\begin{aligned} H(n) &= \binom{n+c_1}{c_1} + \binom{n+c_2-1}{c_2} + \cdots + \binom{n+c_{n+1}-n}{c_{n+1}} \\ &= \binom{n+c_1}{n} + \binom{n+c_2-1}{n-1} + \cdots + \binom{n+c_n-(n-1)}{1} + 1 \end{aligned}$$

and so by [Rob90], Proposition 4.3 we get

$$\begin{aligned} (H(n)_n)_+^+ &= \binom{n+1+c_1}{n+1} + \binom{n+1+c_2-1}{n} + \cdots + \binom{n+1+c_n-(n-1)}{2} + c_n + 1 + 1 \\ &\geq \binom{n+1+c_1}{n+1} + \binom{n+1+c_2-1}{n} + \cdots + \binom{n+1+c_n-(n-1)}{2} + c_{n+1} + 1 + 1 \\ &= \binom{n+1+c_1}{n+1} + \binom{n+1+c_2-1}{n} + \cdots + \binom{n+1+c_n-(n-1)}{2} + \binom{n+1+c_{n+1}-n}{c_{n+1}} + 1 \\ &\geq H(n+1). \end{aligned}$$

□

**Theorem 4.4 (Characterization of Hilbert polynomials)** *Let $P \in \mathbb{Q}[X;\mathbb{N}]$ and $R_0$ an equicharacteristic Artinian local ring. Then the following conditions are equivalent:*

*(i) There exists a standard $R_0$-algebra $S$ such that $h_S = P$,*

*(ii) $P$ admits a Gotzmann development.*

**Proof.** We have just seen that (i) implies (ii) in Corollary 3.7. Assume that $P$ admits a Gotzmann development. Using Corollary 2.10, it is enough to find a function $H : \mathbb{N} \to \mathbb{N}$ such that $H(0) = \lambda_{R_0}(R_0)$, $H(n+1) \leq (H(n)_n)_+^+$ and $H(n) = P(n)$ for $n \gg 0$. This can be made in many ways, but maybe the simplest is to take $H = G[c_1, \ldots, c_s]$ for $n > 0$ and $H(0) = \lambda_{R_0}(R_0)$. □

**Remark 4.5** *The above criterion is effective. For instance, we see that the polynomial $P(X) = 2X$ is not the Hilbert polynomial of any standard $R_0$-algebra with $R_0$ an Artinian local equicharacteristic ring.*

Our aim is now to formulate a more precise version of the characterization Theorem 4.4 which takes into account the number of variables.

**Lemma 4.6** *Let $P \in \mathbb{Q}[X;\mathbb{N}]$ be a polynomial and $c = \deg(P)$. For all $n \geq 0$ consider the Euclidean division*

$$P(n) = \binom{n+c}{c}\gamma(n) + \Gamma(n), \quad 0 \leq \Gamma(n) < \binom{n+c}{c}.$$

*Then $\gamma(n+1) = \gamma(n)$ for all $n \gg 0$ and so $\Gamma$ is asymptotically polynomical.*

**Proof.** Perform the polynomial division

$$P(X) = \binom{X+c}{c}a + Q(X) \quad 0 \leq \deg(Q) < c.$$

Notice that, since $P \in \mathbb{Q}[X;\mathbb{N}]$, $a \in \mathbb{N}$ and $Q \in \mathbb{Q}[X;\mathbb{N}]$. There are two cases to consider:

(a) Assume that the leading coefficient of the polinomial $Q(X)$ is positive. Since the degree of $Q(X)$ is strictly smaller than $c$, we have $0 \leq Q(n) < \binom{n+c}{c}$ for all $n \gg 0$. Hence, taking $\gamma = a$ and $\Gamma(X) = Q(X)$,

$$P(n) = \binom{n+c}{c}\gamma + \Gamma(n)$$

is the Euclidean division for $n \gg 0$, and then $\gamma(n) = \gamma$ for all $n \gg 0$.

(b) If the leading coefficient of the polinomial $Q(X)$ is negative, define then $\Gamma(X) = \binom{X+c}{c} + Q(X)$ and $\gamma = a - 1 \geq 0$. For all $n \gg 0$ we have $Q(n) < 0$ and, since the degree of $Q(X)$ is strictly smaller than $c$, $Q(n) > -\binom{n+c}{c}$. Hence $0 \leq \Gamma(n) < \binom{n+c}{c}$ and therefore

$$P(n) = \binom{n+c}{c}\gamma + \Gamma(n)$$

is the Euclidean division for $n \gg 0$. □

**Definition 4.7** *We will denote by $\gamma(P)$ the limit of the sequence $\{\gamma(n)\}_{n\in\mathbb{N}}$ and by $\Gamma_P \in \mathbb{Q}[X;\mathbb{N}]$ the polynomial associated to the function $\Gamma$. In other words, we have*

$$P(X) = \gamma(P)\binom{X+c}{c} + \Gamma_P(X)$$

*and $P(n) = \gamma(P)\binom{n+c}{c} + \Gamma_P(n)$ is the Euclidean division for all $n \gg 0$. Notice that the polynomical expression of $P(X)$ in terms of $\gamma(P)$ and $\Gamma_P$ is not necessarily the polynomial division, see case (b) in the proof of Lemma 4.6.*

**Example 4.8** *The polynomial $P(X) = X$ has $\gamma(P) = 0$ and $\Gamma_P(X) = X$.*

**Lemma 4.9** *Let $P \in \mathbb{Q}[X;\mathbb{N}]$ and assume that $\Gamma_P$ admits a Gotzmann development. Then $P$ admits a Gotzmann development. Furthermore, if $\gamma(P) \neq 0$ then $\deg(\Gamma_P) < \deg(P)$.*

**Proof.** We may assume that $\gamma(P) \neq 0$ since otherwise the claim is obvious. Let $c = \deg(P)$ and pick any Artinian local equicharacteristic ring $R_0$ such that $\lambda_{R_0}(R_0) = \gamma(P) + 1$. Let $s$ be the length of the Gotzmann development of $\Gamma_P$ and consider $H : \mathbb{N} \to \mathbb{N}$ defined by

$$H(n) = \begin{cases} (\gamma(P)+1)\binom{n+c}{c} & \text{if } n \leq s-1 \\ P(n) & \text{if } n \geq s \end{cases}$$

This function obviously verifies the conditions in Theorem 2.8 (ii) and hence there exists a homogeneous ideal $I \subseteq R_0[X_1,\ldots,X_{c+1}]$ such that $H = H_{R/I}$; it follows that $P = h_{R/I}$ and so it admits a Gotzmann development. Furthemore, by Theorem 3.5 there exist integers $c > c'_1 \geq \ldots \geq c'_p \geq 0$ and $0 \leq q \leq \gamma(P) + 1$ such that

$$P(n) = q\binom{n+c}{c} + \binom{n+c'_1}{c'_1} + \binom{n+c'_2-1}{c'_2} + \cdots + \binom{n+c'_p-(p-1)}{c'_p}$$

for all $n \gg 0$; notice that this is an Euclidean division and hence we must have $q = \gamma(P)$ and

$$\Gamma(X) = \binom{X+c'_1}{c'_1} + \binom{X+c'_2-1}{c'_2} + \cdots + \binom{X+c'_p-(p-1)}{c'_p}$$

that is, $\deg(\Gamma) = c'_1 < c$. □

**Remark 4.10** *The converse does not hold in Lemma 4.9: a counterexample is the polynomial $P(X) = X^2 + 5X - 5$. It admits a Gotzmann development with $s_3 = 2, s_2 = 5$ and $s_1 = 5$, e.g.*

$$P(X) = \binom{X+2}{2} + \binom{X+1}{2} + \binom{X-1}{1} + \binom{X-2}{1} + \binom{X-3}{1}.$$

*Furthermore $P(X) = 2\binom{X+2}{2} + 2X - 7$, that is $\gamma(P) = 2 \neq 0$ and $\Gamma_P(X) = 2X - 7$. But $\Gamma_P$ does not admit a Gotzmann development.*

The following theorem decides whether $P$ is $b$-admissible or not in terms of combinatorical properties of $P$. It will be the main tool used in Section 5 to compute the minimal $b$ for which a function $H$ is $b$-admissible.

**Theorem 4.11** Let $(R_0, \mathfrak{m})$ be an Artinian equicharacteristic local ring, $P \in \mathbb{Q}[X; \mathbb{N}]$ and $b \geq 1$ an integer. Then we have:

(i) If $P$ is $b$-admissible then $b \geq \deg(P) + 1$ and $P$ admits a Gotzmann development.

(ii) If $b \geq \deg(P) + 2$ and $P$ admits a Gotzmann development then $P$ is $b$-admissible.

(iii) $P$ is $(\deg(P) + 1)$-admissible if and only if either of the following conditions holds:

   (a) $0 < \gamma(P) < \lambda_{R_0}(R_0)$ and $\Gamma_P$ admits a Gotzmann development,
   
   (b) $\gamma(P) = \lambda_{R_0}(R_0)$ and $\Gamma_P = 0$.

**Proof.** Let $r = \lambda_{R_0}(R_0)$ and $c = \deg(P)$.

(i) Assume that $P$ is $b$-admissible. By Corollary 3.5 we get that $P$ admits a Gotzmann development. Moreover

$$c = \deg(P) = \dim(R/I) - 1 \leq \dim(R) - 1 = b - 1,$$

that is $b \geq c + 1$.

(ii) Assume now that $P$ admits a Gotzmann development and $b \geq \deg(P) + 2$. We may assume $b = c + 2$. From Corollary 2.11 we deduce that it is enough to show that there exists an admissible function $H : \mathbb{N} \to \mathbb{N}$ such that $H(1) = c + 2$ and $H(n) = P(n)$ for all $n \gg 0$. Since this has been shown in Lemma 4.3, (ii) is proved.

To prove (iii), assume in the first place that $P$ is $(c+1)$-admissible. Then by Theorem 3.5, there exist integers $c > c'_1 \geq \ldots \geq c'_p \geq 0$ and $0 \leq q \leq r$ such that

$$P(n) = q\binom{n+c}{c} + \binom{n+c'_1}{c'_1} + \binom{n+c'_2-1}{c'_2} + \cdots + \binom{n+c'_p-(p-1)}{c'_p}$$

for all $n \gg 0$. Since $c = \deg(P)$ it follows that $q > 0$. Then, this is the Euclidean division of $P(n)$ by $\binom{n+c}{c}$, so we necessarily have $q = \gamma(P)$ and

$$\Gamma_P(X) = \binom{X+c'_1}{c'_1} + \binom{X+c'_2-1}{c'_2} + \cdots + \binom{X+c'_p-(p-1)}{c'_p}.$$

Hence $0 < \gamma(P) \leq r$ and $\Gamma_P$ admits a Gotzmann development. Furthermore if $\gamma(P) = r$, then for $n \gg 0$ we have

$$r\binom{n+c}{c} + \Gamma_P(n) = P(n) = H_{R/I}(n) \leq r\binom{n+c}{c}.$$

Hence $\Gamma_P(n) = 0$ for $n \gg 0$ and so $\Gamma_P = 0$.

Reciprocally, assume that $P$ verifies the conditions in (a) or (b) and let $s$ be the length of the Gotzmann development of $\Gamma_P$ in case (a), $s = 0$ in case (b). To show that $P$ is $(c+1)$-admissible it suffices to constuct a function $H : \mathbb{N} \to \mathbb{N}$ verifying the conditions in Theorem 2.8 (ii) and such that $H(n) = P(n)$ for all $n \gg 0$. It is now immediate that

$$H(n) = \begin{cases} \lambda_{R_0}(R_0)\binom{n+c}{c} & \text{if } n \leq s-1 \\ P(n) & \text{if } n \geq s \end{cases}$$

verifies the conditions required.  $\square$

**Remark 4.12** *Gotzmann developments can be very complicated; for example, the Gotzmann development of $8\binom{n+3}{3}$ has 161427 terms (apply the formulas in Proposition 3.11 with $e_0 = 8$ and $e_1 = e_2 = e_3 = 0$). Notice that if we have an Artinian ring $R_0$ with $\lambda_{R_0}(R_0) = 8$ and $R = R_0[X_1, X_2, X_3, X_4]$, then these expressions are the ones obtained for $h_R$ in parts (i) of Theorem 3.5 and Corollary 3.7 respectively. See also examples 5.6 and 5.7.*

In fact, the expression of $h_{R/I}$ given in Theorem 3.5 is always better than the one obtained in Corollary 3.7:

**Proposition 4.13** *Let $R = R_0[X_1, \ldots, X_b]$, $I \subseteq R_+$ a homogeneous ideal and*

$$h_{R/I}(X) = q\binom{X+b-1}{b-1} + \binom{X+c_1'}{c_1'} + \binom{X+c_2'-1}{c_2'} + \cdots + \binom{X+c_p'-(p-1)}{c_p'}$$

$$= \binom{X+c_1}{c_1} + \binom{X+c_2-1}{c_2} + \cdots + \binom{X+c_s-(s-1)}{c_s}$$

*the expressions of $h_{R/I}$ obtained in Theorem 3.5 and Corollary 3.7 respectively. Define $p_i$ and $s_i$ as in these two results; then for all $i \geq 1$ it holds $p_i \leq s_i$.*

**Proof.** We may assume $q \neq 0$, hence $c_1 = b - 1$. We have

$$q\binom{n+b-1}{b-1} + G[c_1', \ldots, c_p'](n) = G[c_1, \ldots, c_s](n)$$

for all $n \gg 0$. The result is obvious for $b = 1$. Notice that for $b \geq 2$ we have

$$\Delta G[c_1, \ldots, c_s](n) = G[c_1 - 1, \ldots, c_{s_2} - 1](n)$$

and

$$\Delta \left( q\binom{n+b-1}{b-1} + G[c_1', \ldots, c_p'](n) \right) = q\binom{n+b-2}{b-2} + G[c_1' - 1, \ldots, c_{p_2}' - 1](n),$$

hence by induction on $b$ it is enough to show that $p \leq s$.

We will proceed by induction on $p$. If $p = 0$ the statement is clear. Assume $p > 0$ and let us distinguish two cases:

*Case (1):* $c_p' = 0$. Then for $n \gg 0$ we have

$$q\binom{n+b-1}{b-1} + G[c_1', \ldots, c_p'](n) = q\binom{n+b-1}{b-1} + G[c_1', \ldots, c_{p-1}'](n) + 1$$

$$= G[b_1, \ldots, b_t](n) + 1 \quad \text{by Lemma 4.9}$$

$$= G[b_1, \ldots, b_t, 0](n).$$

So we must have $s = t + 1$. Since by induction hypothesis we have $p - 1 \leq t$, we get $p \leq s$.

*Case (2):* $c_p' > 0$. Let $t = c_p' < b - 1 = c_1$. Then for $n \gg 0$ we have

$$\Delta^t \left( q\binom{n+b-1}{b-1} + G[c_1', \ldots, c_p'](n) \right) = \Delta^t(G[c_1, \ldots, c_s](n))$$

that is

$$q\binom{n+b-1-t}{b-1-t} + G[c_1' - t, \ldots, c_p' - t](n) = G[c_1 - t, \ldots, c_{s_{t+1}} - t](n)$$

therefore by case (1) we have $p \leq s_{t+1} \leq s$. □

**Proposition 4.14 (Sharp lower bound for the Hilbert function)** *Let $(R_0, \mathfrak{m})$ be an Artinian equicharacteristic local ring and $S$ a standard $R_0$-algebra. Let $c_1, \ldots c_s$ be the Gotzmann coefficients of $h_S$; then we have*

$$H_S(n) \geq G[c_1, \ldots, c_s](n)$$

*for all $n \geq 1$. Besides, the bound above is sharp.*

**Proof.** The proof uses exactly the same argument given in the proof of Theorem 3.5 to show that $H^1_{R_+}(R/I)_n = 0$ for all $n \geq p - 1$. The sharpness is consequence of Lemma 4.3. □

Notice that in the 1-dimensional case the bound above reduces to the well-known result

$$H(n) \geq \begin{cases} n+1 & \text{if } n \leq e_0(S) - 1 \\ e_0(S) & \text{if } n \geq e_0(S) - 1. \end{cases}$$

**Remark 4.15** *In the local case, when $n = 1$, since $c_1 = dim A - 1$ and $s_d = e_0(A)$, Proposition 4.14 reduces to Krull's height theorem: If $(A, \mathfrak{m})$ is not a regular ring (hence $e_0(A) \geq 2$), then the minimal number of generators of $\mathfrak{m}$, i.e., $H^0_A(1)$, is $\geq dim A + 1$.*

By induction on $t$ we can also prove, using Corollary 3.10 in [Sta78]

**Proposition 4.16** *Let $(R_0, \mathfrak{m})$ be an Artinian equicharacteristic local ring, $S$ a standard $R_0$-algebra and assume that $depth(S) \geq t$. Let $c_1, \ldots c_s$ be the Gotzmann coefficients of $h_S$; then we have for all $0 \leq i \leq t$*

$$\Delta^i H_S(n) \geq G[c_1 - i, \ldots, c_{s_{i+1}} - i](n)$$

*for all $n \geq 1$. Besides, the bound above is sharp.*

## 5 Admissibility of functions. Ideals with a given Hilbert function

Let $H : \mathbb{N} \to \mathbb{N}$ be an asymptotically polynomical function and $(R_0, \mathfrak{m})$ any Artinian equicharacteristic local ring. Our aim in this section is to give an algorithm to decide whether $H$ is an admissible function. It seems that the natural way to encode $H$ should be to give a finite number of values of $H$, say $H(0), H(1), \ldots, H(n_0)$, and a polynomial $h(X) \in \mathbb{Q}[X; \mathbb{N}]$ such that $H(n) = h(n)$ for all $n > n_0$. From Theorem 2.10 we know that $H$ is admissible if and only if it verifies the conditions in (ii); the problem is to verify these conditions in a finite number of steps. The theory of Gotzmann developments will provide us with a method to do so. Furthermore, in the case that $H$ is admissible we compute the minimal value $b$ for which $H$ is $b$-admissible. We also describe how to get a generating system for an ideal $I \subseteq R_0[X_1, \ldots, X_b]$ such that $H = H_{R/I}$; in the case $R_0 = k$ it will be a minimal generating system.

Let us begin by giving an algorithm to decide whether a polynomial $P \in \mathbb{Q}[X]$ is an admissible polynomial.

## 5.1 Algorithm to compute the Gotzmann development

Here we give an algorithm to compute, if it exists, the Gotzmann development of a polynomial. The strategy is to compute first the normalized Hilbert-Samuel coefficients and then compute from them the Gotzmann coefficients using Proposition 3.11. The following proposition provides a triangular equation system in $e_1, \ldots, e_c$ and a criterion to decide whether $P \in \mathbb{Q}[X; \mathbb{N}]$.

**Proposition 5.1** *Let $h \in \mathbb{Q}[X]$ be a polynomial of degree $c$ and $e_0, \ldots, e_c$ its normalized Hilbert-Samuel coefficients. Then we have for all $0 \leq i \leq c - 1$.*

$$\sum_{k=0}^{i} \binom{i}{k} e_{c-k} = (-1)^c h(-i-1)$$

*In particular, $h \in \mathbb{Q}[X; \mathbb{N}]$ if and only if its leading coefficient is a positive multiple of $1/c!$ and $h(-1), \ldots, h(-c) \in \mathbb{Z}$.*

**Proof.** We have $h(X) = \sum_{k=0}^{c} (-1)^{c-k} e_{c-k} \binom{X+k}{k}$, and notice that for all $0 \leq i \leq c - 1$

$$\binom{k-i-1}{k} = \begin{cases} 0 & \text{if } i < k \\ (-1)^k \binom{i}{k} & \text{if } i \geq k \end{cases}$$

So we get for all $0 \leq i \leq c - 1$

$$h(-i-1) = \sum_{k=0}^{c} (-1)^{c-k} e_{c-k} \binom{k-i-1}{k} = \sum_{k=0}^{i} (-1)^c e_{c-k} \binom{i}{k}$$

as we wanted to prove. The second part is consequence of the fact that the system matrix is unipotent upper triangular with coefficients in $\mathbb{N}$. □

**ALGORITHM GOTZTST:**

**INPUT:** A polynomial $P \in \mathbb{Q}[X]$.
**OUTPUT:** The normalized Hilbert-Samuel coefficients of $P$; the Gotzmann development of $P$, in case it exists, and its length $s$.

**Step 1:** Make sure that the leading coefficient of $P$ is positive and compute the values $P(-1), \ldots, P(-c)$, where $c = \deg(P)$ together with $e_0 = c! \cdot$(leading coefficient of $P$). If any of them is not an integer then by Proposition 5.1 $P \notin \mathbb{Q}[X; \mathbb{N}]$ and so it can not admit a Gotzmann development, hence we stop here. Assume all of them are integer.

**Step 2:** Solve the triangular system of equations of Proposition 5.1 in order to obtain $e_c, e_{c-1}, \ldots, e_1$. Notice that the coefficients of this system are the entries in the Pascal triangle, so the system will be computationally well-behaved. We can obtain the coefficients of every equation by shifting the preceding ones to the left and to the right and adding.

**Step 3**: Once we have computed the Hilbert-Samuel coefficients, solve the triangular system of equations of Proposition 3.11 to obtain $s_1, \ldots, s_c$. The only fact we must check at each step is whether $s_i \leq s_{i+1}$, for all $c - 1 \geq i \geq 1$. If this holds, then $P$ admits a Gotzmann development and so it is an admissible polynomial. Let $s = s_1$: then the Gotzmann coefficients are obtained as $c_1 = \cdots = c_{s_{c+1}} = c$, $c_{s_{c+1}+1} = \cdots = c_{s_c} = c - 1, \ldots, c_{s_2+1} = \cdots = c_s = 0$.

Notice that an alternative way to compute $e_0$ should be to use the equation $e_0 - e_1 + \cdots + (-1)^{c-1}e_{c-1} + (-1)^c e_c = P(0)$ in Step 2. Here we should check that $P(0) \in \mathbb{Z}$ and the leading coefficient of $P$ is positive in Step 1 instead of computing $e_0 = $ (leadcoeff)$\cdot c!$. The advantage is that we do not need to compute $c!$ with this method.

## 5.2 Algorithm to check $b$-admissibility

Fix an Artinian equicharacteristic local ring $(R_0, \mathfrak{m})$ and a set of data $(r, i_1, \ldots, i_{n_0}, h(X))$ describing a function $H : \mathbb{N} \to \mathbb{N}$. We assume that $r = \lambda_{R_0}(R_0)$, $i_n \in \mathbb{N}$ for all $1 \leq i \leq n_0$ and $h(X) \in \mathbb{Q}[X]$. Hence we have $H(0) = r$, $H(n) = i_n$ for all $1 \leq n \leq n_0$ and $H(n) = h(n)$ for $n \geq n_0 + 1$. Notice that any admissible function can be encoded in this way.

We are going to use Remark 2.12, Corollary 2.11, Theorem 4.11, Proposition 5.1, Proposition 4.6, Lemma 4.9, Theorem 4.4 and Proposition 3.9 in order to check whether $(r, i_1, \ldots, i_{n_0}, h(X))$ describe an admissible function and, in such case, to compute the minimal $b$ for which $H$ is admissible.

**ALGORITHM $b$-ADM:**

**INPUT:** $r = \lambda_{R_0}(R_0)$ and an asymptotically polynomical function $H : \mathbb{N} \to \mathbb{N}$, encoded as $H = (i_1, \ldots, i_{n_0}, h(X))$, with $n_0 \geq 0$, $i_j \in \mathbb{N}$ and $h(X) \in \mathbb{Q}[X; \mathbb{N}]$.
**OUTPUT:** A decision about whether $H$ is admissible or not. If $H$ is admissible, then the minimal $b$ for which $H$ is $b$-admissible is computed.

**Step 1:** Let $c = \deg(h)$. If $H(1) < c + 1$ then by Remark 2.12 and Theorem 4.11 (i) $H$ is not admissible, in this case stop. Hence we have $H(1) \geq c + 1$.

**Step 2:** Perform Step 1 of algorithm GOTZTST in order to determine if $h \in \mathbb{Q}[X; \mathbb{N}]$: make sure that the leading coefficient of $h$ is positive and compute the values $h(-1), \ldots, h(-c)$ together with $e_0 = c! \cdot$(leading coefficient of $h$). If any of them is not an integer then by Proposition 5.1 $h \notin \mathbb{Q}[X; \mathbb{N}]$ and hence $H$ is not admissible, so we stop here.

If all of them are integer, compute the normalized Hilbert-Samuel coefficients of $h$, $e_0, \ldots, e_c$ as in Step 2 of algorithm GOTZTST.

**Step 3:** If $H$ were $b$-admissible, by Remark 2.12 we would have $H(1) \leq br$ and by Theorem 4.11 $b \geq c+1$. Hence, set $b_{min} = \max\{]H(1)/r[, c+1\}$, where $]x[$ denotes the least integer greater or equal than $x$. If $b_{min} > c + 1$ then skip to Step 6, since Theorem 4.11 (iii) does not apply here. Otherwise we want to study the $(c+1)$-admissibility of $H$, hence we are going to apply Theorem 4.11 (iii). Recall that in the proof of Proposition 4.6 we have $Q(X) = h(X) - e_0\binom{X+c}{c}$. If $H$ is $(c+1)$-admissible $\Gamma_h$ must admit a Gotzmann development, hence by Lemma 4.9 together with the construction of $\gamma(h)$ and $\Gamma_h$ in the proof of Proposition 4.6 we must have $\gamma(h) = e_0$ and $\Gamma_h = Q$.

**Step 4:** If $e_0 > r$ then by Theorem 4.11 (iii) we skip to Step 6.

**Step 5:** Now we are going to check whether $\Gamma_h = Q$.

   **5.1** If $e_1 = \ldots = e_c = 0$, i.e. $Q = 0$, set $p = 0$ and go to Step 7.

- **5.2** If the first nonzero $e_i$ verifies $(-1)^i e_i < 0$, i.e. the leading coefficient of $Q$ is negative, this means that $\Gamma_h \neq Q$, by the proof of Proposition 4.6. Therefore, set $b_{min} = b_{min} + 1$ and skip to Step 6.

- **5.3** Assume that the first nonzero $e_i$ verifies $(-1)^i e_i > 0$, i.e. the leading coefficient of $Q$ is positive. If $e_0 = r$ set $b_{min} = b_{min} + 1$ and skip to Step 6. If $e_0 < r$ compute the Gotzmann development of $\Gamma_h = Q$ as in algorithm GOTZTST. Notice that the normalized Hilbert-Samuel coefficients of $\Gamma_h$ are $e'_j = (-1)^i e_{i-j}$, $0 \leq j \leq c - i$, so we do not need to perform Steps 1 and 2 in algorithm GOTZTST.

  If the Gotzmann development of $\Gamma_h$ does not exist, set $b_{min} = b_{min} + 1$ and skip to Step 6. If it exists, let $p$ be its length and skip to Step 7.

**Step 6:** Compute the Gotzmann development of $h$. Notice that the normalized Hilbert-Samuel coefficients have already been computed in Step 2, so we can go directly to Step 3 in algorithm GOTZTST. If $h$ does not admit a Gotzmann development then stop here: by Theorem 4.4 $H$ is not admissible. Otherwise let $p$ be its length.

**Step 7:** Make sure that $i(H) = n_0 + 1$. For this, compare $i_{n_0}$ and $h(n_0)$. If they coincide delete the superfluous data $i_{n_0}$ and iterate until $i_{n_0} \neq h(n_0)$.

**Step 8:** Define $m = \max\{i(H), p\}$ and for $n$ between $0$ and $m$ compute the Euclidean division $H(n) = \binom{n+b_{min}-1}{b_{min}-1} q(n) + r(n)$, see Proposition 3.9. Check at every step that $(q(n), r(n)) \leq (q(n-1), (r(n-1)_{n-1})^+_+)$. If this holds for every $1 \leq i \leq m$ then $H$ is an admissible function and the minimal value of $b$ such that $H$ is admissible is $b = b_{min}$. Furthermore $H$ will be $b$-admissible for all $b \geq b_{min}$. If this fails then

- **8.1** If $b_{min} = c + 1$, set $b_{min} = b_{min} + 1$ and skip to Step 6. Notice that this time it will not be necessary to perform Step 7 since $i(H)$ is already computed.

- **8.2** If $b_{min} = H(1)$ then thanks to Corollary 2.11 $H$ is not admissible, so we stop here.

- **8.3** Otherwise set $b_{min} = b_{min} + 1$ and skip to the head of Step 8.

### 5.3 Algorithm to constuct an ideal with a given Hilbert function

Given a $b$-admissible function $H = (r, i_1, \ldots, i_{n_0}, h(X))$ and a composition series $\mathcal{J} = \{0 = J_0 \subseteq J_1 \subseteq \ldots \subseteq J_r = R_0\}$ in $R_0$ we will show how to construct the ideal $I_{H,\mathcal{J}} \subseteq R_0[X_1, \ldots, X_b]_+$ which appeared in the proof of Theorem 2.8.

**Proposition 5.2** *Let $k$ be a field, $I \subseteq k[X_1, \ldots, X_b] = R$ a segment ideal and $H = H_{R/I}$. Then it holds:*

*(i) For all $n \geq 1$, $H(n+1) = (H(n)_n)^+_+$ if and only if $I_{n+1} = R_1 I_n$*

*(ii) $R_1 I_n$ is generated as a $k$-vector space by the last $\binom{n+1+b-1}{b-1} - (H(n)_n)^+_+$ monomials in $R_n$.*

**Proof.** (i) follows from [BH93], Proposition 4.2.8 and Macaulay's theorem. (ii) follows from [BH93], Lemma 4.2.5 and Proposition 2.4. □

**ALGORITHM EFEC:**

**INPUT:** $r = \lambda_{R_0}(R_0)$, a composition series $\mathcal{J} = \{0 = J_0 \subseteq J_1 \subseteq \ldots \subseteq J_r = R_0\}$, a function $H$ encoded as in algorithm $b$-ADM and $b \geq 1$ an integer for which we know that $H$ is admissible.
**OUTPUT:** The $\mathcal{J}$-segment ideal $I_{H,\mathcal{J}} \subseteq R_0[X_1, \ldots, X_b]_+$.

**Step 1:** Compute $i(H)$ as in Step 7 of algorithm $b$-ADM.

**Step 2:** Set $p =$ length of the Gotzmann development of $\Gamma_h$ if $b = \deg(h) + 1$ and $p =$ length of the Gotzmann development of $h$ if $b > \deg(h) + 1$.

**Step 3:** Let $m = \max\{i(H), p\}$. By Proposition 3.9 we know that $I_{H,\mathcal{J}}$ is generated in degrees at most $m$. For $0 \leq n \leq m$ let us compute the Euclidean division $H(n) = Nq(n) + r(n)$, where $N = \binom{n+b-1}{b-1}$, and the values $(g_1(n), g_2(n)) = (q(n-1) - q(n), (r(n-1)_{n-1})^+_+ - r(n)) \geq (0,0)$ for $1 \leq n \leq m$. Then, if we set

$$\nu(n) = \min\{r(n), (r(n-1)_{n-1})^+_+\} = \begin{cases} r(n) & \text{if } g_2(n) \geq 0 \\ (r(n-1)_{n-1})^+_+ & \text{if } g_2(n) < 0, \end{cases}$$

the generators of $I_{H,\mathcal{J}}$ in degree $n$ are

$$J_{r-q(n)-1}X^{\lambda_1}, \ldots, J_{r-q(n)-1}X^{\lambda_{r(n)}}, J_{r-q(n)}X^{\lambda_{r(n)+1}}, \ldots, J_{r-q(n)}X^{\lambda_N} \quad \text{if } g_1(n) > 1,$$

$$J_{r-q(n)-1}X^{\lambda_1}, \ldots, J_{r-q(n)-1}X^{\lambda_{\nu(n)}}, J_{r-q(n)}X^{\lambda_{r(n)+1}}, \ldots, J_{r-q(n)}X^{\lambda_N} \quad \text{if } g_1(n) = 1,$$

$$J_{r-q(n)}X^{\lambda_{r(n)+1}}, \ldots, J_{r-q(n)}X^{\lambda_{r(n)+g_2(n)}} \quad \text{if } g_1(n) = 0,$$

Notice that the ones we have skipped among the generators which appeared in the proof of Theorem 2.8 are superfluous by Proposition 5.2. Also $r(n) + g_2(n) = (r(n-1)_{n-1})^+_+$.

Let us finally make some remarks about the case where $R_0$ is a field:

**Remark 5.3** *If $k$ is a field, $H$ is an admissible function and $b = H(1)$, then the generating system obtained above is*

$$\bigcup_{n=2}^{m} \{X^{\lambda_{H(n)+1}}, \ldots, X^{\lambda_{(H(n-1)_{n-1})^+_+}}\}$$

*and it is a minimal generating system for $I_H$.*

*Applying then [ERV91], Corollary 2.7, we can compute for every homogeneous ideal $J \subseteq R = k[X_1, \ldots, X_b]$ a bound for its minimal number of generators that will depend only on $H = H_{R/J}$. Namely, if $m = \max\{i(H), s\}$, where $s$ is the length of the Gotzmann development of $h$, then*

$$\nu(J) \leq \sum_{n=1}^{m} (H(n) - (H(n-1)_{n-1})^+_+).$$

*For example, if $J$ is any ideal having as Hilbert function the one given in Example 5.5 below then $\nu(J) \leq 4$.*

In this case we can also compute the zero-th local cohomology group of $R/I_H$: let $I \subseteq R = k[X_1, \ldots, X_b]$ be a homogeneous ideal and $J \subseteq R$ be the homogeneous ideal such that $H^0_{R_+}(R/I) = J/I$. We have $J_n = I_n$ for all $n > a_0(R/I)$. In other words, $J = I^{sat}$ is the saturation of $I$: it is the biggest homogeneous ideal containing $I$ and having the same Hilbert polynomial and verifies $\text{depth}(R/J) \geq 1$.

**Lemma 5.4** *Let $s$ be the length of the Gotzmann development of $h_{R/I}$ and $m = \max\{i(R/I), s\}$. Then we have:*

  (i) *If $I$ is a monomial ideal then $J$ is monomial too.*

  (ii) *$J$ is generated in degrees at most $s$.*

  (iii) *If $I$ is a segment ideal and $\alpha \in R$ is a monomial, then for all $n \geq 1$, $\alpha R_n \subseteq I$ if and only if $\alpha X_1^n \in I$.*

  (iv) *If $I$ is a segment ideal then $J = \bigcup_{n \geq 1}(I : X_1^n)$.*

  (v) *If $I$ is a segment ideal, then $J$ is generated by $\bigcup_{n=0}^{s}\bigcup_{i=0}^{m-n}(I_{n+i} : X_1^i)$.*

**Proof.** (i) Let $\alpha \in J$ be an element which we may assume to be homogeneous. Write $\alpha = t_1 + \cdots + t_r$ as a sum of terms; then it is enough to show that every $t_i \in J$.

Since $\alpha \in J$ there exists $n \in \mathbb{N}$ such that $\alpha X^\lambda \in I$ for all multi-indices $\lambda$ with $|\lambda| = n$. That is to say $t_1 X^\lambda + \cdots + t_r X^\lambda \in I$ for all $|\lambda| = n$. Notice that this is a sum of terms and $I$ is a monomial ideal, hence $t_i X^\lambda \in I$ for all $|\lambda| = n$ and for all $i$, i.e., $t_i \in J$ for all $i$. (ii) is Proposition 3.9. (iii) follows from the fact that if $\beta, \gamma \in R_n$ are monomials with $\beta > \gamma$, then $\alpha\beta > \alpha\gamma$. Since $X_1^n$ is the biggest monomial in $R_n$, we are done. (iv) and (v) follow from the first three claims. □

**Example 5.5** *Let $R_0 = k$, $k$ any field. We will use algorithm b-ADM to prove that $H = (1, 4, 10, 19, X^2 + 3X)$ is admissible and to compute the minimal value of $b$ for which $H$ is admissible. Since $H(1)/\lambda_{R_0}(R_0) \leq b_{min} \leq H(1)$ we know that we will get $b_{min} = 4$. Then for $b = 4$ we will compute a minimal generating system of $I_H$ using algorithm EFEC, as well as the zero-th local cohomology group of $R/I_H$.*

**Step 1:** $c = \deg(h) = 2$ and $H(1) = 4 \geq c + 1$.

**Step 2:** $h(-1) = -2, h(-2) = -2$, the leading coefficient of $h$ is positive and $e_0 = 2! \cdot 1 = 2$: $h \in \mathbb{Q}[X; \mathbb{N}]$.
$e_2 = h(-1) = -2$
$e_1 + e_2 = h(-2) = -2$, so $e_1 = 0$.

**Step 3:** $b_{min} = \max\{H(1), c+1\} = 4$, so we skip to Step 6.

**Step 6:** Let us compute the Gotzmann development of $h$ using algorithm GOTZTST:
$s_3 = e_0 = 2$, $s_2 = \binom{s_3+1}{2} - e_1 = 3 \geq s_3$, $s_1 = -\binom{s_3+1}{3} + \binom{s_2+1}{2} + e_2 = 3 \geq s_2$.
Therefore $h$ is an admissible polynomial: it has Gotzmann development

$$h(X) = \binom{X+2}{2} + \binom{X+1}{2} + \binom{X-1}{1}$$

and $p = s = 3$.

**Step 7:** Since $H(3) = 19 \neq 18 = h(3)$ we have $i(H) = n_0 + 1 = 4$.

**Step 8:** $m = \max\{i(H), s\} = 4$.

$(q(0), r(0)) = (1, 0)$
$(q(1), r(1)) = (1, 0) \leq (q(0), (r(0)_0)_+^+) = (1, 0)$
$(q(2), r(2)) = (1, 0) \leq (q(1), (r(1)_1)_+^+) = (1, 0)$
$(q(3), r(3)) = (0, 19) \leq (q(2), (r(2)_2)_+^+) = (1, 0)$
$(q(4), r(4)) = (0, 28) \leq (q(3), (r(3)_3)_+^+) = (0, 31)$
Hence $H$ is an admissible function and $b_{min} = 4$.

Let us compute $I_H \subseteq k[X_1, X_2, X_3, X_4]$.
$(g_1(1), g_2(1)) = (0, 0)$, hence $I_1 = 0$.
$(g_1(2), g_2(2)) = (0, 0)$, hence $I_2 = 0$.
$(g_1(3), g_2(3)) = (1, -19)$ and $\nu = 0$, hence

$$I_3 = \langle J_{1-q(3)} X^{\lambda_{r(3)+1}}, \ldots, J_{1-q(3)} X^{\lambda_N} \rangle = \langle X^{\lambda_{20}} \rangle,$$

since $J_0 = k$ and $N = \binom{3+3}{3} = 20$.
The ordered monomials of degree 3 are
$X_1^3 > X_1^2 X_2 > X_1 X_2^2 > X_2^3 > X_1^2 X_3 > X_1 X_2 X_3 > X_2^2 X_3 > X_1 X_3^2 > X_2 X_3^2 > X_3^3 > X_1^2 X_4 > X_1 X_2 X_4 > X_2^2 X_4 > X_1 X_3 X_4 > X_2 X_3 X_4 > X_3^2 X_4 > X_1 X_4^2 > X_2 X_4^2 > X_3 X_4^2 > X_4^3.$
We must delete the first 19 and keep the remaining one; that is, we have to keep $X_4^3$.

$(g_1(4), g_2(4)) = (0, 3)$, hence

$$I_4 = \langle J_{1-q(4)} X^{\lambda_{r(4)+1}}, \ldots, J_{1-q(4)} X^{r(4)+g_2(4)} \rangle = \langle X^{\lambda_{29}}, X^{\lambda_{30}}, X^{\lambda_{31}} \rangle,$$

since $J_0 = k$.
The ordered monomials of degree 4 are
$X_1^4 > X_1^3 X_2 > X_1^2 X_2^2 > X_1 X_2^3 > X_2^4 > X_1^3 X_3 > X_1^2 X_2 X_3 > X_1 X_2^2 X_3 > X_2^3 X_3 > X_1^2 X_3^2 > X_1 X_2 X_3^2 > X_2^2 X_3^2 > X_1 X_3^3 > X_2 X_3^3 > X_3^4 > X_1^3 X_4 > X_1^2 X_2 X_4 > X_1 X_2^2 X_4 > X_2^3 X_4 > X_1^2 X_3 X_4 > X_1 X_2 X_3 X_4 > X_2^2 X_3 X_4 > X_1 X_3^2 X_4 > X_2 X_3^2 X_4 > X_3^3 X_4 > X_1^2 X_4^2 > X_1 X_2 X_4^2 > X_2^2 X_4^2 > X_1 X_3 X_4^2 > X_2 X_3 X_4^2 > X_3^2 X_4^2 > X_1 X_4^3 > X_2 X_4^3 > X_3 X_4^3 > X_4^4.$
We must delete the first 28 and keep the first three of the remaining ones; that is, we have to keep $X_1 X_3 X_4^2, X_2 X_3 X_4^2$ and $X_3^2 X_4^2$.
Since $I_H$ is generated in degrees $n \leq 4$, the minimal generating system for $I_H$ is:

$$I_H = (X_4^3, X_1 X_3 X_4^2, X_2 X_3 X_4^2, X_3^2 X_4^2).$$

The system of generators for the ideal $J$ such that $H^0_{R_+}(R/I) = J/I$ is

$$J = (X_4^3, X_3 X_4^2).$$

Notice that in the case $R_0 = k$ it is useless to compute the ideal $I_H \subseteq k[X_1, \ldots, X_b]$ for values of $b > b_{min}$; we would get the same ideal with the addition of the $b - b_{min}$ last variables.

**Example 5.6** Let $R_0 = k[T]/(T^4) = k[\varepsilon]$, $k$ any field, and let $H = (4X + 4)$. Applying the preceding algorithms we get that $H$ is an admissible function, $h(X) = 4\binom{X+1}{1}$ in terms of the Hilbert-Samuel coefficients so $\gamma(h) = 4 = \lambda_{R_0}(R_0)$ and $\Gamma_h = 0$; since $m = max\{i(H), p\} = 0$ we get $b_{min} = 2$.

*Fix $\mathcal{J} = \{0 \subseteq (\varepsilon^3) \subseteq (\varepsilon^2) \subseteq (\varepsilon) \subseteq R_0\}$ as a composition series in $R_0$.*

*For $b = 2$ we obtain $I_{H,\mathcal{J}} = 0$.*

*For $b = 3$ we obtain $I_{H,\mathcal{J}} = (\varepsilon^3 X_1, \varepsilon^3 X_2, \varepsilon^2 X_3, \varepsilon^2 X_1^2, \varepsilon^2 X_1 X_2, \varepsilon^2 X_2^2, \varepsilon^2 X_2 X_3, \varepsilon X_1 X_3^2, \varepsilon X_2 X_3^2, \varepsilon X_3^3, \varepsilon X_1^3 X_3, \varepsilon X_1^2 X_2 X_3, \varepsilon X_1^2 X_2^3, \varepsilon X_1 X_2^4, \varepsilon X_2^5, \varepsilon X_1^6, \varepsilon X_1^5 X_2, \varepsilon X_1^4 X_2^2, X_2^2 X_3^5, X_1 X_3^6, X_2 X_3^6, X_3^7, X_1^2 X_2 X_3^5, X_1^4 X_3^5, X_2^6 X_3^4).$*

*If we wanted to obtain an ideal generated only by monomials without coefficients in $\mathfrak{m}$ we would have to perform the algorithm with $b > H(1)$.*

**Example 5.7** *Let $R_0 = k[T]/(T^5) = k[\varepsilon]$, $k$ any field, and let $H = (9, 11, 3X + 5)$. Applying the preceding algorithms we get that $H$ is an admissible function, $h(X) = 3\binom{X+1}{1} + 2$ in terms of the Hilbert-Samuel coefficients, so $\gamma(h) = 3$ and $\Gamma_h = 2$, $m = max\{i(H), p\} = 2$ and $b_{min} = 2$.*

*Fix $\mathcal{J} = \{0 \subseteq (\varepsilon^4) \subseteq (\varepsilon^3) \subseteq (\varepsilon^2) \subseteq (\varepsilon) \subseteq R_0\}$ as a composition series in $R_0$.*

*For $b = 2$ we obtain $I_{H,\mathcal{J}} = (\varepsilon^4 X_2, \varepsilon^4 X_1^2, \varepsilon^3 X_2^2)$.*

*For $b = 3$ we obtain $I_{H,\mathcal{J}} = (\varepsilon^3 X_1, \varepsilon^3 X_2, \varepsilon^3 X_3, \varepsilon^2 X_1^2, \varepsilon^2 X_1 X_2, \varepsilon^2 X_2^2, \varepsilon^2 X_1 X_3, \varepsilon^2 X_2 X_3, \varepsilon X_3^2, \varepsilon X_1^2 X_3, \varepsilon X_1 X_2 X_3, \varepsilon X_2^2 X_3, \varepsilon X_1^2 X_2^2, \varepsilon X_1 X_2^3, \varepsilon X_2^4, \varepsilon X_1^3 X_3, \varepsilon X_1^2 X_2 X_3, \varepsilon X_1^5, \varepsilon X_1^4 X_2, X_3^5, X_1 X_2 X_3^4, X_2^2 X_3^4, X_1^3 X_3^4, X_2^5 X_3^2).$*

# References

[BH93]   W. Bruns and J. Herzog. *Cohen-Macaulay rings*, volume 39 of *Cambridge studies in advanced mathematics*. Cambridge University Press, 1993.

[CLO92]  D. Cox, J. Little, and D. O'Shea. *Ideals, varieties, and algorithms*. Undergraduate texts in Mathematics. Springer-Verlag, 1992.

[ERV91]  J. Elias, L. Robbiano, and G. Valla. Number of generators of ideals. *Nagoya Math. J.*, 123:39–76, 1991.

[Got78]  G. Gotzmann. Eine Bedingung für die Flachheit und das Hilbertpolynom eines graduierten Ringes. *Math. Z.*, 158:61–70, 1978.

[Gre89]  M. Green. Restrictions of linear series to hyperplanes, and some results of Macaulay and Gotzmann. In E. Ballico and C. Cilliberto, editors, *Algebraic Curves and Projective Geometry Proceedings, Trento (1988)*, volume 1389 of *Lecture Notes in Mathematics*, pages 76–86. Springer-Verlag, 1989.

[Hil90]  D. Hilbert. Über die Theorie der algebraischen Formen. *Math. Ann.*, 36:473–534, 1890.

[HIO88]  M. Herrmann, S. Ikeda, and U. Orbanz. *Equimultiplicity and blowing up*. Springer-Verlag, 1988.

[Hoa93a] L. T. Hoa. Reduction numbers of equimultiple ideals. (preprint), 1993.

[Hoa93b]  L. T. Hoa. Two notes on coefficients of the Hilbert-Samuel polynomial. (preprint), 1993.

[Mac16]  F. S. Macaulay. *The algebraic theory of modular systems*. Cambridge University, 1916.

[Mac27]  F. S. Macaulay. Some properties of enumeration in the theory of modular systems. *Proc. London Math. Soc.*, 26:531–555, 1927.

[Mar93]  T. Marley. The reduction number of an ideal and the local cohomology of the associated graded ring. *Proc. Amer. Math. Soc.*, 117:335–341, 1993.

[Mat86]  H. Matsumura. *Commutative ring theory*, volume 8 of *Cambridge studies in advanced mathematics*. Cambridge University Press, 1986.

[Mum66]  D. Mumford. *Lectures on curves on an algebraic surface*, volume 59 of *Annals of Mathematics Studies*. Princeton University Press, 1966.

[Rob90]  L. Robbiano. Introduction to the theory of Hilbert functions. In *The Curves Seminar at Queen's, Vol VII*, number 85 in Queen's Papers in Pure and Appl. Math. Queen's Univ., Kingston, Canada, 1990.

[Sam51]  P. Samuel. La notion de multiplicité en algèbre et en géometrie algébrique. *J. de Math.*, XXX:159–274, 1951.

[Ser65]  J.P. Serre. *Algèbre locale. Multiplicités*, volume 11 of *Lecture Notes in Mathematics*. Springer-Verlag, 1965.

[Sta78]  R. P. Stanley. Hilbert functions of graded algebras. *Adv. in Math.*, 28:57–83, 1978.